\documentclass[10pt]{article}
\usepackage{times,amssymb,euscript,epsfig,latexsym,amsmath}
\usepackage[normalem]{ulem}
\pdfoutput=1
\usepackage{epsfig,graphicx}
\usepackage{color}
\begin{document}
\newcommand{\todo}[1]{\vspace{5 mm}\par \noindent
\framebox{\begin{minipage}[c]{0.95 \textwidth}
\tt #1 \end{minipage}}\vspace{5 mm}\par}
\newtheorem{theorem}{Theorem}[section]
\newtheorem{corollary}{Corollary}
\newtheorem{lemma}[theorem]{Lemma}
\newtheorem{proposition}{Proposition}
\newtheorem{proof}{Proof}[section]
\date{}
\title{ \Large{On the Dynamics of a Nutation Ring Damper}}
\author{Adrian Luna \\
Banavara N. Shashikanth\footnote{Mechanical and Aerospace Engineering
Department, MSC 3450, PO Box 30001, New Mexico State University,
Las Cruces, NM 88003, USA. E-mail:shashi@nmsu.edu}}

\maketitle
\pagenumbering{gobble}
Nutation triggered by external disturbances is an undesirable phenomenon for a spinning satellite. Passive decay in nutation angle must be achieved in such cases and benefits missions involving spin-stabilized satellites or mass-constrained satellites like CubeSats. Motivated by previous work of other authors, the dynamics of a  general nutating satellite with an attached ring containing a viscous fluid is analyzed in a theoretical mechanics framework. The fluid in the ring is modeled as a rigid point mass and is then extended to a distributed rigid slug. The nonlinear equations of motion of the system, under the assumption of zero external torques, are studied in a body-fixed frame. Equilibrium points and system trajectories are studied on the angular momentum sphere. The settling time of the dynamical system after an initial disturbance is investigated. The energy dissipation is modeled using an approximate drag coefficient. Fluid dynamic effects are further examined through CFD simulations in the Navier-Stokes setting using the open-source software package ANSYS Fluent. Applications of the analysis and simulations to CubeSats are discussed.
\newpage
\pagenumbering{arabic}
\section{Introduction.}

Nutation ring dampers have been used by the space industry to dampen the wobble of the spin axis of a satellite on encountering a disturbance. The dampers are constructed using small masses or viscous fluids to dissipate unwanted rotational oscillations. The dampening is achieved passively which avoids the need for an active stabilization that can require control algorithms and powered devices. 

   Attitude control is essential for the success of satellite missions for maintaining desired orientation in irregular environments.  Spin stabilized spacecraft, utilizing the gyroscopic effect of the rotating  spacecraft, were favored during initial satellite development due to their innate ability to maintain the spacecraft in a particular direction.  Although it is still the more commonly used method for attitude stabilization of satellites, 
the addition of a nutation damper is necessary for the spin vector to maintain alignment with the spacecraft’s desired spin axis \cite{mobley1971attitude}. Rotating spacecraft in the absence of external forces may not experience precession but nutation can still occur as it is caused internally. Spin stabilized spacecraft can also experience nutation wobble due to deviations in the principal axis \cite{roldugin2023wobble}. Modern large satellites now utilize three-axis stabilization; however, many small modular satellites would benefit from the reliability and passive means for attitude control.

Liquid ring nutation dampers in spacecraft were first used in the Pioneer 1 lunar probe in 1958 \cite{herzl1971introduction}. The dampers work by transferring the transverse angular momentum to the primary spin axis which causes the satellite to achieve single axis rotation. Koval and Bhuta analyzed the dynamics of the damper by investigating the energy dissipation through their fluid model and experimental validation \cite{bhuta1966viscous}. They conclude that the damper's light weight allows for the device to not interfere with the dynamics of the satellite. Alfriend modeled the liquid as a rigid slug which was found to agree with experimental results \cite{alfriend1974partially}. Alfriend and Spencer proposed parameters maximizing the performance of nutation dampers in their investigation of nutation decay expressions \cite{alfriend1983comparison}. The ring damper can be placed perpendicular or parallel to the spin axis, without a significant effect on the performance. The perpendicular displacement from the center of mass plane however is an important parameter regarding performance \cite{alfriend1974partially}. The dynamics of the damper can be split into large and small angle damping modes. The damper in low nutation angles can cause the fluid to spread out uniformly from its slug-like behavior which can separate and behave like point masses \cite{cartwright1963circular}.

   Concurrently, there has been a growth in the nanosatellite industry \cite{swartwout2013first,villela2019towards}. Nanosatellites are miniaturized satellites that can fill excess payloads and test instruments. This has led to the broad use of CubeSats which are a class of nanosatellites that have been used for a variety of missions for education and research \cite{heidt2000cubesat}. The development of CubeSats has impacted the accessibility of satellite research due to low cost and modular construction. Their size and flexibility make them useful for short term research or tests before being implemented into large satellites \cite{chin2008cubesat,poghosyan2017cubesat}. These satellites are typically operated in low-Earth orbits where environmental disturbances can be significant. Spacecraft in such orbits can be subject to aerodynamic drag and radiation pressure effects that can cause undesired surface forces \cite{klinkrad1998orbit}. CubeSats are measured in standard unit modules defined by a cubic volume of side length of 10 cm and mass of 1-2 kilograms \cite{heidt2000cubesat}. Modular satellites have the ability to dock and interact during missions to offset their low volumes. However, the small size of each unit restricts power and propulsion capacity which creates the need for energy efficient control methods. Three-axis stabilization methods allow for fixed orientations, but the attitude devices often used are rotational devices as alternative options such as thrusters are dependent on fuel. Rotational devices such as reaction wheels store angular momentum which can experience wobble in the event of a disturbance \cite{dennehy2019survey}. Use of these devices can cause unwanted nutation that can negatively impact the performance of instruments and the stability of the spacecraft. The deployment of passive stabilization has the potential to increase the spacecraft’s lifetime and enhance mission capabilities.

 Papers on theoretical modeling and analysis of nutation ring dampers have been few and far between [5,6,15-20]
%
This paper investigates the theoretical framework and implementation of a nutation ring damper, with potential applications for nanosatellites.  It highlights certain important features such as the angular momentum sphere seemingly missed in previous work. The analysis is accompanied by numerical simulations to study the damper’s performance during expected operating conditions. The damper model studied consists of a hollow ring with a fluid slug inside the ring. Three different models of the slug are examined: (a) Point mass model, (b) Distributed fluid model and (c) CFD simulations
Most of the paper focuses on the dynamic behavior in models (a) and (b). A novel feature of this paper is the presentation of Navier-Stokes simulations of the liquid inside the ring using the open-source CFD software Ansys Fluent. Modern computational power was not available to the authors of previous work. Applications to CubeSats are discussed using typical CubeSat parameters.   
      
The organization of the paper is as follows. In Section 2, the configuration and the dynamics of the point mass model is described, and the expressions for the system momenta and kinetic energy in a body-fixed frame are presented. The phase space of the system is discussed and MATLAB simulations of trajectories on the angular momentum sphere and time plots are presented. In Section 3, the distributed rigid slug is considered. The equations of motion remain the same, but the inertia tensors are modified corresponding to the distributed mass. MATLAB simulations are presented. In Section 4, CFD simulations of the Navier-Stokes equations for a viscous fluid slug are presented using the open-source program ANSYS Fluent.  Finally, the paper is summarized with a discussion of future research directions.

\section{Point Mass Slug.} 

\subsection{Dynamical framework and theoretical analysis.}


  The configuration of the damper is that of a circular ring, which contains a slug, that encircles the satellite. Assume the satellite is modeled as a rigid cylinder of radius $R$, which is also the radius of the ring, and length $L$. The composite cylinder-ring-slug body will be referred to as the `ring-slug' system since from the dynamics perspective the presence of the cylinder is accounted for simply by including its inertia coefficients. Similarly, in most cases, the word `ring' will imply ring+cylinder. 

    As shown in Fig. \ref{rs}, let $XYZ$ denote a stationary frame centered at $O$, and let $xyz$ denote a body-fixed frame with its $z$-axis parallel to the cylinder's axis.  The $xyz$ frame with unit vectors $(e_1,e_2,e_3)$ is centered at $o$, the center of mass of the ring-slug system,  assumed to lie in a cross-sectional plane of the cylinder. The $e_1,e_2$ axes lie in the plane of the ring and the $e_1$ axis is such that it is always directed towards the slug. The slug is allowed to move freely and frictionlessly, relative to the ring walls, at an angular velocity $\dot{\beta}e_3$ about an axis perpendicular to the plane of the ring. Let the entire system's center of mass lie at point $\mathcal{O}$ at a distance $d$ from $o$ along the $x$-axis. The axis of rotation passes through $\mathcal{O}$.
The presence of the slug causes the center of mass of the ring to move from its center $\mathcal{O}$ to $o$ to, with 
\begin{align}
d&=\frac{m_s}{m_s+m_r}R  \label{eq:cm1} \\
\Rightarrow m_rd&=m_s(R-d) \label{eq:cm2}
\end{align}
Since in most satellite applications, $m_s << m_r$, the above shows that $d$ is small. However, because it scales with the mass ratio $m_s/m_r$, it is important in this theoretical framework not to set $d=0$ since this would then correspond to the nonphysical cases of $m_s=0$ or $m_r=\infty$. The effects of small $d$ can be obtained  only after a careful order of magnitude analysis. Since the inclusion of $d$ does not cause any added complexity in the analysis, it will be retained throughout.

It is assumed that there  are no external forces or torques acting on the combined ring-slug system. The former implies that the center of mass of the ring-slug system moves with a constant velocity which may be assumed to be zero. Since there is no linear velocity, the value of $b$ in Fig. \ref{rs} can be set to zero and points $O$ and $o$ taken as coincident.

\subsubsection{Body-fixed frame: Momenta and Kinetic Energy.}
The system angular momentum and kinetic energy expressions referred to the body-fixed frame are
\begin{align}
\mathbf{h}(t)&=\mathbf{I_r} \mathbf{\Omega} + \mathbf{I_s} \mathbf{\Omega}_s, \quad K.E.=\frac{1}{2}(\mathbf{\Omega})^T \mathbf{I_r} \mathbf{\Omega}+\frac{1}{2}(\mathbf{\Omega}_s)^T  \mathbf{I_s}  \mathbf{\Omega}_s \label{eq:appambf}
\end{align}
where $\mathbf{\Omega}$ and $\mathbf{\Omega}_s$ are the ring and slug angular velocities, respectively, written in a body-fixed frame and related as
\begin{align}
\mathbf{\Omega}_s&=\mathbf{\Omega}+\dot{\beta}e_3. \label{eq:slugavbf}
\end{align}
$\mathbf{I_r}$ and $\mathbf{I_s}$ are the moments of inertia matrices of the ring and slug, respectively, about the center of mass of the system. As mentioned above, the $\mathbf{I_r}$ entries need to take into account the presence of the cylinder (representing the satellite). The final expressions for these matrices are:
\begin{align*}
\bf{I_r}&=\begin{bmatrix}
I_{xr} & 0 & 0\\
0 & I_{yr}& 0\\
0 & 0 & I_{zr}
\end{bmatrix}=
\begin{bmatrix}
I_{xr,0} & 0 & 0\\
0 & I_{yr,0} +m_r d^2& 0\\
0 & 0 & I_{zr,0} +m_rd^2
\end{bmatrix},
\end{align*}
where
\begin{align*}
I_{xr,0}&=\frac{m_r R^2}{4}+\frac{m_r L^2}{12}, \quad 
I_{yr,0}=\frac{m_r R^2}{4}+\frac{m_r L^2}{12} \quad
I_{zr,0}=\frac{m_r R^2}{2}
\end{align*}
and 
\begin{align}
\bf{I_s}=
\begin{bmatrix}
I_{xs} & 0 & 0\\
0 & I_{ys} & 0\\
0 & 0 & I_{zs} 
\end{bmatrix}=\begin{bmatrix}
0 & 0 & 0\\
0 & m_s(R-d)^2 & 0\\
0 & 0 & m_s(R-d)^2 
\end{bmatrix} \label{eq:pmslugI}
\end{align}
%
\noindent{\it Remark.} Expressions (\ref{eq:appambf}) are the same as the ones used by previous authors. But it is worth noting that they are obtained without imposing the conservation of linear momentum. There is a small non-zero amount of time-varying linear momentum which is ignored.  In the exactly conserved linear momentum model, the system rotates like a dumbbell with a single angular velocity $\mathbf{\Omega}_s$.


\subsubsection{Equations of motion and Phase Space.} In the absence of any external torques acting on the combined system, the total system angular momentum about the point $o$ in the stationary frame  is conserved, i.e.
\begin{align}
\frac{d \bar{\mathbf{h}}}{dt}&=0
\end{align}
\begin{figure}[h]
\centering
\includegraphics[width=2in]{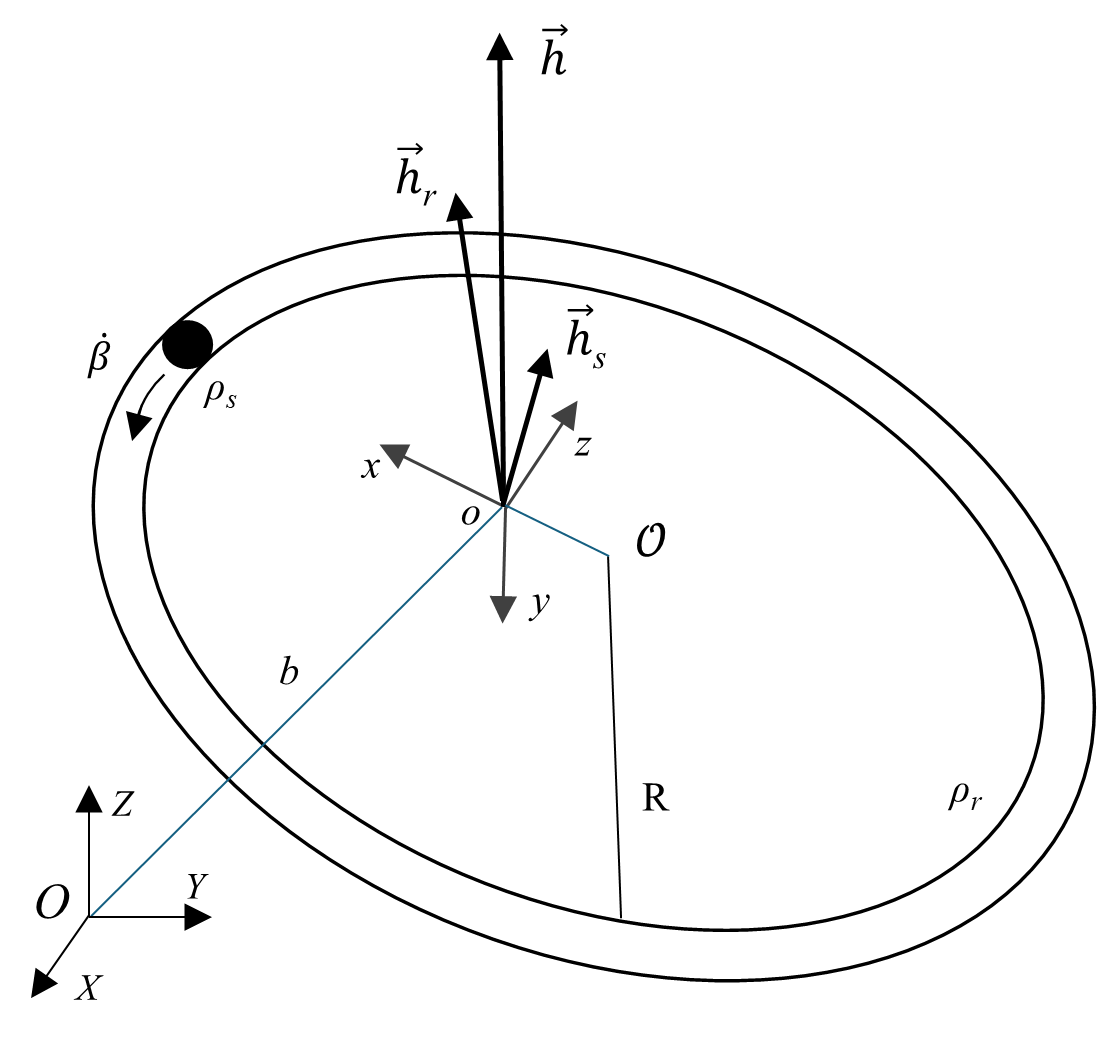}
\caption[Schematic sketch of the ring damper with the spatial angular momentum vectors]{A schematic sketch of the ring damper when the fluid slug is modeled as a point mass (black dot). Note that the satellite body, which is modeled as a cylinder encircled by the ring, is not shown. The spatial angular momentum vectors of the ring, the slug and the total  system, are also shown.}
\label{rs}
\end{figure}
 The evolution equations of $\mathbf{h}$ are then obtained using the standard rules of transformation of vectors from a stationary to a rotating frame. Namely,  $\bar{\mathbf{X}}(t)$ in the stationary frame and $\mathbf{X}(t)$ in the body-fixed frame are related by the relation: $\bar{\mathbf{X}}(t)=\mathbf{A}(t) \mathbf{X}(t), \dot{\bar{\mathbf{X}}}=\mathbf{A}(t) \left(\dot{\mathbf{X}} + \mathbf{\Omega} \times \mathbf{X} \right)$.
where $\mathbf{A}(t) \in \operatorname{SO}(3)$ is the instantaneous rotation matrix that takes frame $XYZ$ to $xyz$ \cite{IntMechMars}. This leads to 
\begin{align}
\frac{d{\bar{\bf h}}}{dt}&=\mathbf{A}(t) \left(\frac{d{\bf h}}{dt}+\mathbf{\Omega_s} \times \bf{h} \right)
\Rightarrow \frac{d{\bf h}}{dt}=-\mathbf{\Omega_s} \times \bf{h}=-\mathbf{\Omega} \times \mathbf{h}-\dot{\beta} e_3 \times \mathbf{h}, \label{eq:angmomeqn}  
\end{align}
%

      Note that the angular momentum equation (\ref{eq:angmomeqn}) can be written in terms of $\mathbf{h}$ alone by inverting (\ref{eq:appambf}) and expressing $\mathbf{\Omega}_s$ in terms of $\mathbf{h}$, using (\ref{eq:slugavbf}),
resulting in
\begin{align}
\frac{d\mathbf{ h}}{dt}&= \mathbf{h} \times \left[(\mathbf{I_r}+\mathbf{I_s})^{-1} (\mathbf{h} +\mathbf{I_r}\dot{\beta} e_3)\right] \label{eq:angmomeqn2}
\end{align}

The form (\ref{eq:angmomeqn}) of the angular momentum equation clearly shows that
\begin{align*}
\frac{d}{dt}( \mathbf{h} \cdot \mathbf{h})&=0,
\end{align*}
using a standard vector identity. This shows that
i.e. {\it the magnitude of the angular momentum $\mathbf{h}(t)$ in the body-fixed frame is conserved by the dynamics, even in the presence of dissipation.} 

  In other words, as the system evolves, the tip of the angular momentum vector always lies on a sphere of radius $\mid \mathbf{h}(0)\mid^2$ in $\mathbb{R}^3$, where $\mathbf{h}(0)$ is the vector of initial values of the angular momentum. This is called the {\it angular momentum sphere} of the system and the quantity $\mid \mathbf{h} \mid^2$ is called a {\it Casimir} function on the phase space of the system  \cite{IntMechMars}. In the absence of dissipation, i.e with $\ddot{\beta}=0$, $\mathbb{R}^3$ is also the phase space of the system which is foliated by these spheres. 
\begin{figure}[h]
\centering
\includegraphics[width=2.5in]{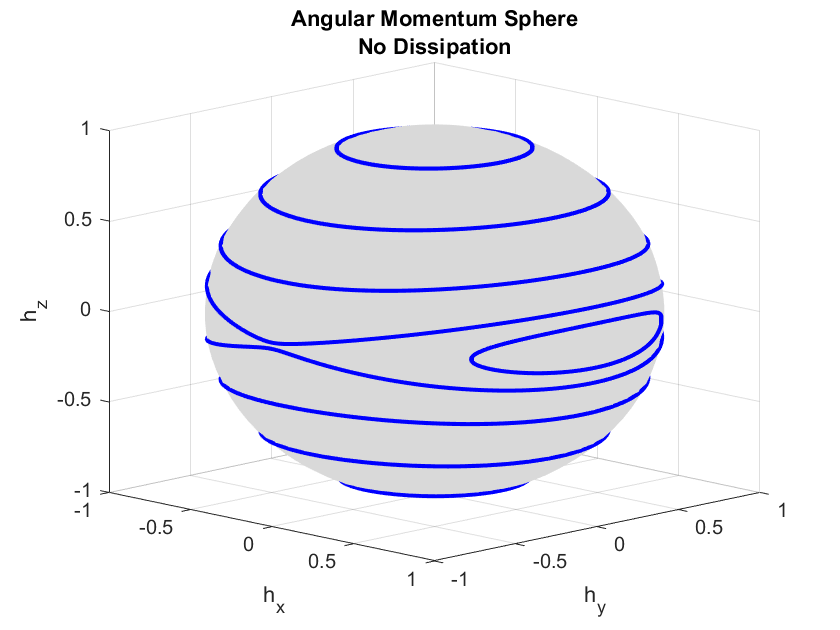}
\caption[Phase space of the system with no dissipation]{Phase trajectories of the system with no dissipation plotted in MATLAB 2022a. Note the well-known features of dissipationless, free rigid body stability \cite{IntMechMars}, namely, rotation about the long axis ($h_z$-axis) is stable, rotation about the intermediate axis ($h_y$-axis) is unstable and rotation about the short axis ($h_x$-axis) is conditionally stable. }
\label{zerodiss}
\end{figure}

 In the presence of dissipation, the phase space is no longer $\mathbb{R}^3$ since now $\dot{\beta}$ is also a dynamical variable of the system whose evolution is coupled with the evolution of $\mathbf{h}(t)$. 
 The phase space can be identified with $\mathbb{R}^4$ but it is important to note that the constraint of the angular momentum sphere still holds, since this is simply a consequence of the form of equation (\ref{eq:angmomeqn}). The trajectories obtained on the sphere in the presence of dissipation are therefore the projections of the actual phase space trajectories in $\mathbb{R}^4$.  To obtain the complete equations of motion of the system, the angular momentum equations (\ref{eq:angmomeqn}) have to be appended with an evolution equation $\dot{\beta}$ which requires a model for the internal frictional forces. 

 \paragraph{Frictional forces and dissipation.} The kinetic energy dissipation is caused by the action and reaction pair of the tangential  frictional forces between the slug and the walls of  the ring.  Note that while friction slows down the slug, the reaction force by the slug on the walls actually results in a slight increase in the angular velocity of the ring. However, the net  result is a decrease in the system kinetic energy. 

     The simplest way to model the tangential frictional force is as a linear function of the slug's linear velocity relative to the ring walls. Introducing a `drag coefficient' $C_d$, write the retarding force of the ring on the slug, in the stationary frame, as 
\begin{align*}
\bar{\mathbf{F}}_{r,s}&=-C_d (R-d) \dot{\beta}(t) \bar{e}_\beta,
\end{align*}
where $\bar{e}_\beta$ is the unit vector in the plane of the ring and tangent to the ring, positive in the direction of $\dot{\beta}$ shown in Fig. \ref{rs}. The negative sign makes the force direction opposite to the slug motion. The corresponding torque about the center of mass $o$ is 
\begin{align*}
\bar{\tau}_{r,s}&=(R-d)\bar{e}_1 \times  \bar{\mathbf{F}}_{r,s}=-C_d (R-d)^2 \dot{\beta}(t) \bar{e}_3,
\end{align*}
where $(\bar{e}_1,\bar{e}_2,\bar{e}_3)$ are the unit vectors $(e_1,e_2,e_3)$ written in the stationary frame. Since the slug's $\bar{e}_3$-component angular momentum about $o$, due to this relative velocity, is $I_{zs} \dot{\beta} \bar{e}_3$, applying
\begin{align*}
\bar{\tau}_{r,s}&=\frac{d}{dt} (I_{zs} \dot{\beta}) \bar{e}_3=I_{zs} \ddot{\beta} \bar{e}_3,
\end{align*}
we get 
\begin{align*}
\ddot{\beta}&=-C_d (R-d)^2 \dot{\beta}/I_{zs}=-C_d \dot{\beta}/m_s
\end{align*}
 Now, Newton's Third Law stipulates that 
\begin{align}
\bar{\mathbf{F}}_{r,s}+\bar{\mathbf{F}}_{s,r}&=0, \label{eq:fsr}
\end{align}
where $\bar{\mathbf{F}}_{s,r}$ is the reaction frictional force exerted by the slug on the ring.  Since the reaction force acts at the same location, it has the same moment arm about $o$ and its torque is equal and opposite, i.e 
\begin{align*}
\bar{\tau}_{r,s}&=-\bar{\tau}_{s,r},
\end{align*}
leading to the relation
\begin{align}
I_{zr}\dot{\Omega}_z&=\bar{\tau}_{s,r}\Rightarrow \dot{\Omega}_z=C_d (R-d)^2 \dot{\beta}/I_{zr} \label{eq:omegazbeta}
\end{align}
Substituting for $h_z$, using (\ref{eq:appambf}),  some fairly straightforward calculations lead to
the final ODE for the $\dot{\beta}$ variable:
\begin{align}
\ddot{\beta}=-\dot{\beta}C_d (R-d)^2\frac{(I_{zs}+I_{zr})}{I_{zs}I_{zr}}+\frac{h_x h_y}{(I_{yr}+I_{ys})I_{zs}}-\frac{h_x h_y}{(I_{xr}+I_{xs})I_{zs}} \label{eq:ddotbeta}
\end{align}
To summarize, the equations of motion of the system are given by (\ref{eq:angmomeqn}) and (\ref{eq:ddotbeta}).

\subsection{Projected Phase Portraits and Time Plots.}
Phase portraits projected on to the angular momentum sphere and time integration of the governing equations of motion are performed utilizing MATLAB 2022a to better understand the behavior of the system. The parameters are scaled to fit a small satellite to replicate conditions CubeSats operate in: Slug mass=0.005 kg, Damper Radius=0.05 m, Friction Coefficient=1.63, Cylinder mass, radius, height=2 kg, 0.05m, 0.05m. The mass of the slug is a tuned parameter which should be within a percent of the satellite's total mass. The initial angular velocity of the system is $\bf{\Omega}=[100,0,400]$, where the disturbance is 25 percent of the satellite's spin rate.

  An examination of equations (\ref{eq:angmomeqn}) and (\ref{eq:ddotbeta}) shows that there are six equilibrium points given by:
\begin{align*}
\mathbf{h_{eq,1}}&=[ 0,0,\pm \mid {\bf h}\mid,0]^T, \quad \mathbf{h_{eq,2}}=[ 0,\pm \mid {\bf h}\mid,0,0]^T, \quad  \mathbf{h_{eq,3}}=[ \pm \mid {\bf h}\mid,0,0,0]^T
\end{align*}
where $\mid {\bf h} \mid=constant$ is the value of the radius of the angular momentum sphere as determined by the initial conditions. Since linear stability analysis is not very insightful, we directly proceed to examine phase portraits projected on to the angular momentum sphere and time histories of various quantities. 
\begin{figure}[h]
\centering
\includegraphics[width=2.5in]{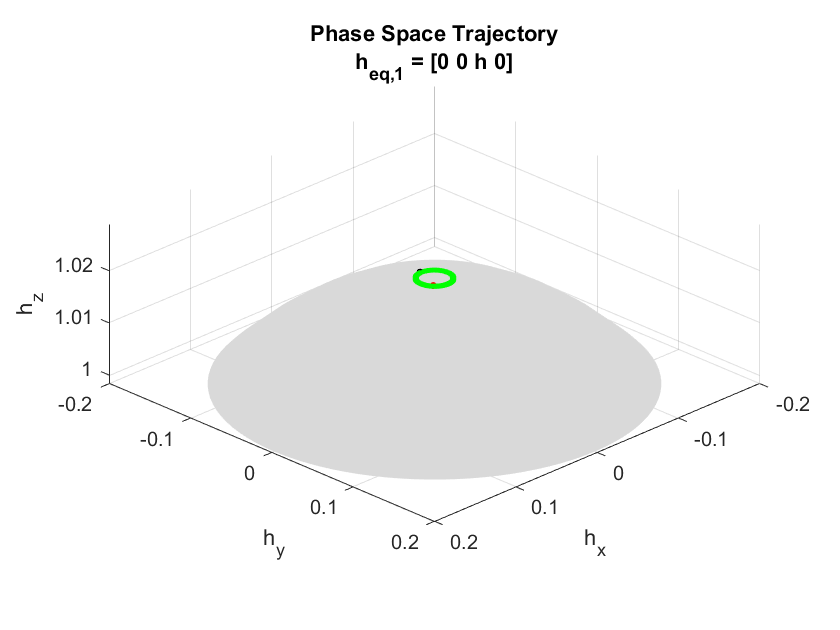}
\includegraphics[width=2.5in]{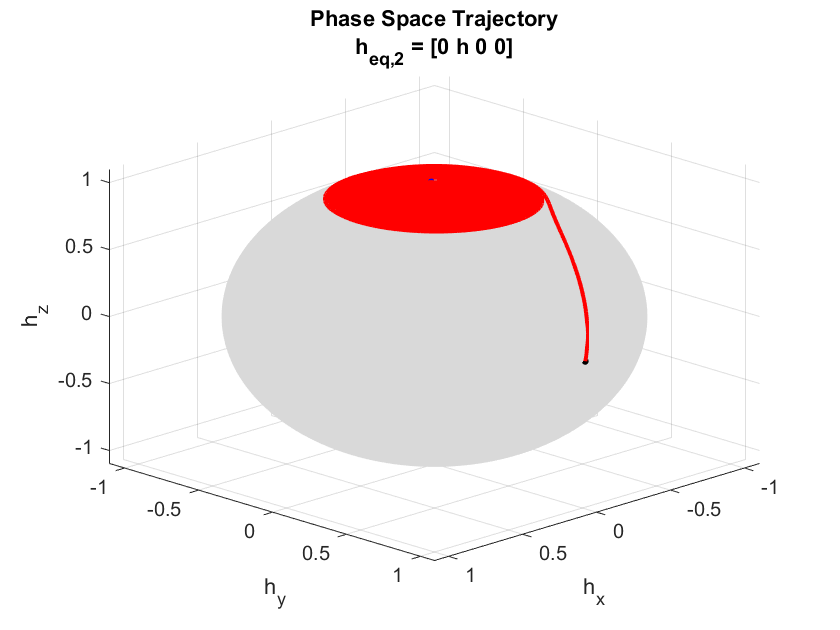}
\includegraphics[width=2.5in]{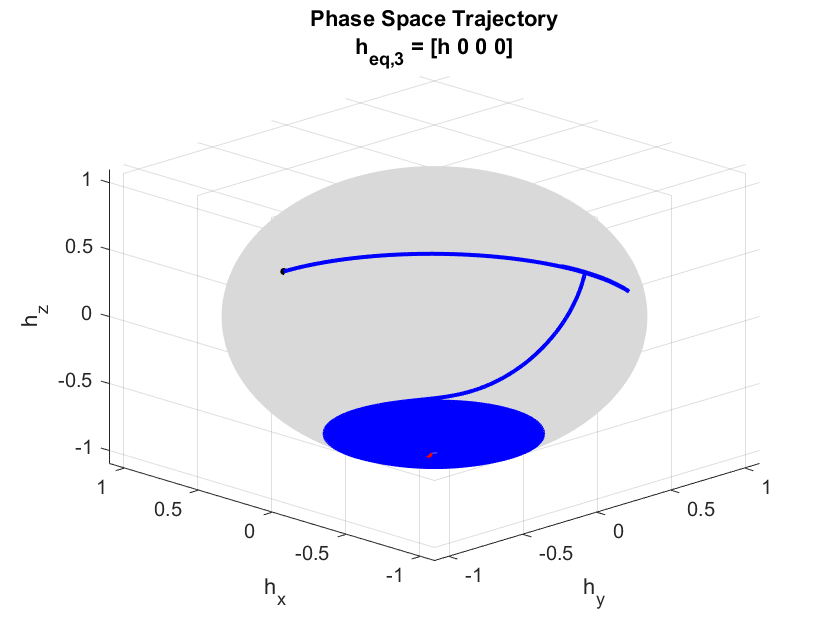}
\includegraphics[width=2.5in]{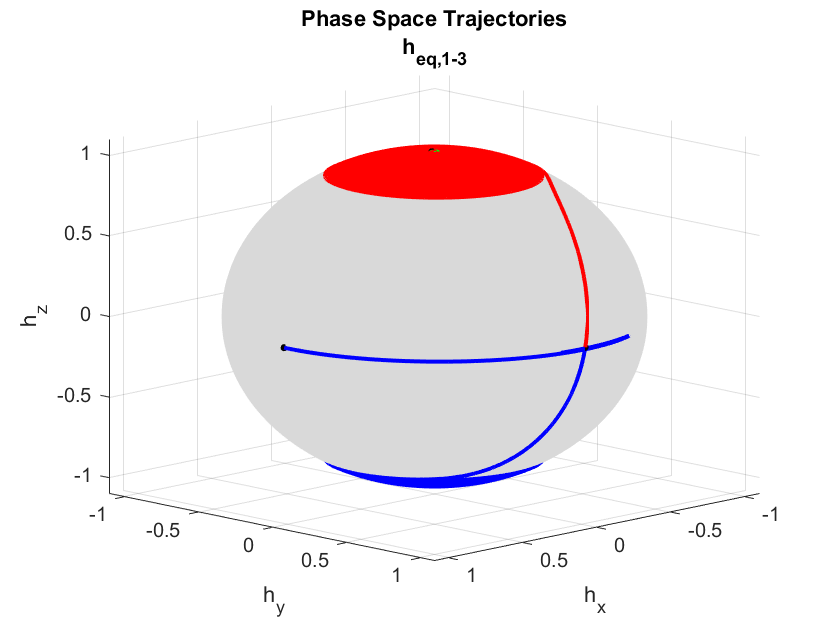}
\caption[Trajectories of the linearized equilibrium points]{The top row figures and the bottom left figure show trajectories on the sphere of the dissipative system starting from close to equilibrium points on each of the three positive axes, the last figure is a composite of these trajectories. }
\label{phasetraj}
\end{figure}
The projected dynamics in the vicinity of the $\mathbf{h_{eq,1}}$ shows stable spiral behavior as seen in Fig. \ref{phasetraj}, top left. This is consistent with the fact that rotation of the structure about this axis has the largest moment of inertia. Interestingly, the trajectory does not continue to move closer to the equilibrium point and so the damper does not fully dampen the transverse angular momentum. 
The trajectory of the second equilibrium point immediately deviates from the initial point. The trajectory moves away from the $h_y$-axis and begins to spiral toward the $h_z$-axis. This particular axis is the intermediate axis of the geometry which would explain its instability. The system rotating about its minor axis also deviates from the equilibrium point. This behavior is different as the system approaches the $h_y$-axis initially. From the $h_y$-axis, the system follows the behavior seen in the second equilibrium point before spiraling towards the $h_z$-axis. These trajectories show that while the damper dissipates the unwanted transverse angular momentum and the kinetic energy, keeping the total magnitude of the angular momentum constant, the overall stability of the system is also influenced by the stability characteristics of the dissipationless, rotating body (Fig. \ref{zerodiss}).

\begin{figure}[h]
\centering
\includegraphics[width=2.4in]{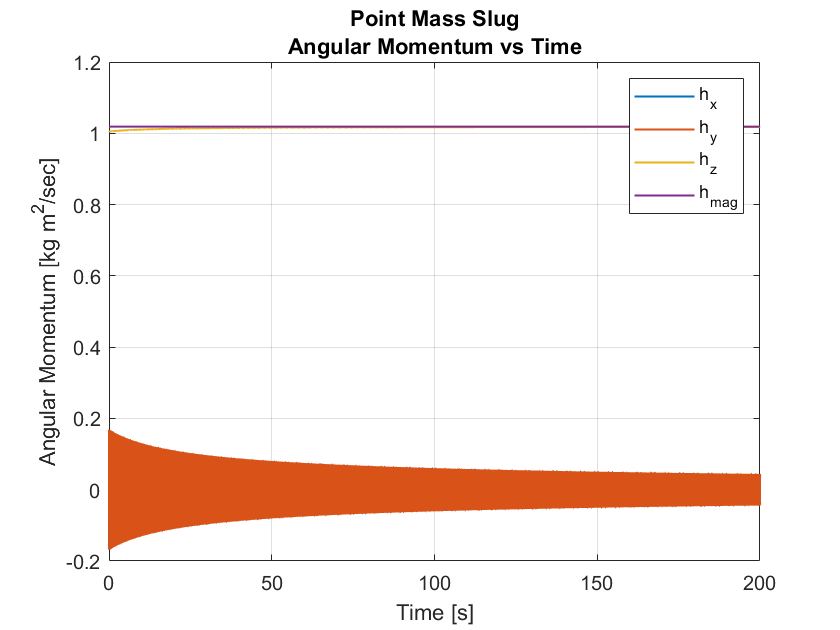}
\includegraphics[width=2.5in]{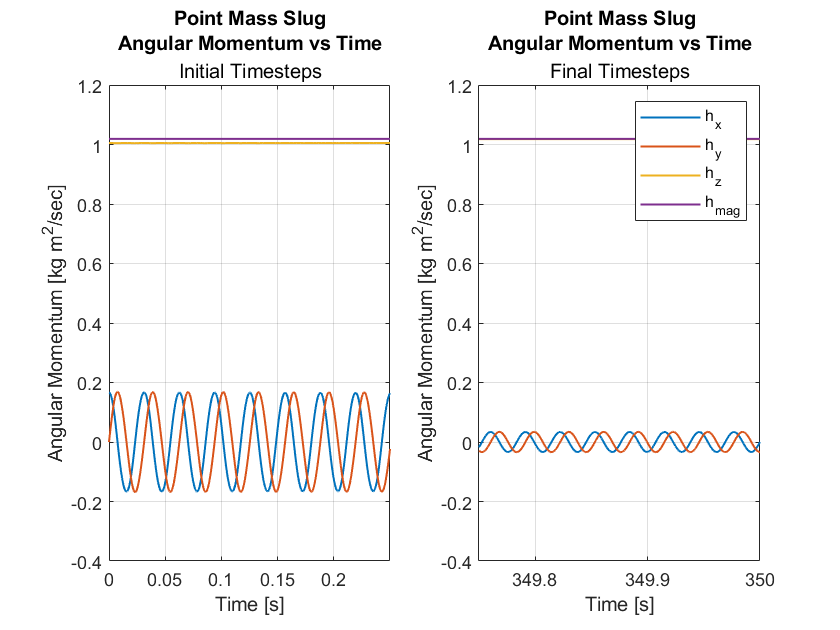}
\caption[Angular momenta components in time (point-mass)]{Plot on the left showing the conservation of the total magnitude of the angular momentum and the variation of each of its components in time. Plots on the right, on much shorter sub-intervals, clearly show the oscillations in the $h_x$ and $h_y$ components.}
\label{angmomplotspm}
\end{figure}
The time evolution of the angular momentum of the system in the body-fixed frame is shown in Fig. \ref{angmomplotspm} (left). The total magnitude of angular momentum is plotted to show its constant value as well as display its relative value to the other components. The transverse angular momentum consists of $x$ and $y$ axis momenta that oscillate throughout the integrated time range. The frequency of the oscillations creates the dense portion of the plot where the amplitude of the oscillations decay as time increases. The oscillations occur due to the wobble motion that is caused by the disturbance. The nutation damper causes the transverse angular momentum to transfer to the polar angular momentum seen in the $z$-component as it increases.
The isolated plots shown in Fig. \ref{angmomplotspm} (right) allow for the comparison of transverse angular momentum during the initial and final time intervals. The reduction of the amplitudes of the transverse angular momentum is significant but does not completely dampen the unwanted angular momentum. The point mass slug offsets the mass distribution from the center of the satellite which causes the wobble to approach a non-zero value. Other works on nutation dampers mention that a distributed slug increases the performance of the damper and reduces the offset of mass distribution \cite{bhuta1966viscous,alfriend1983comparison}.
\begin{figure}[h]
\centering
\includegraphics[width=2in]{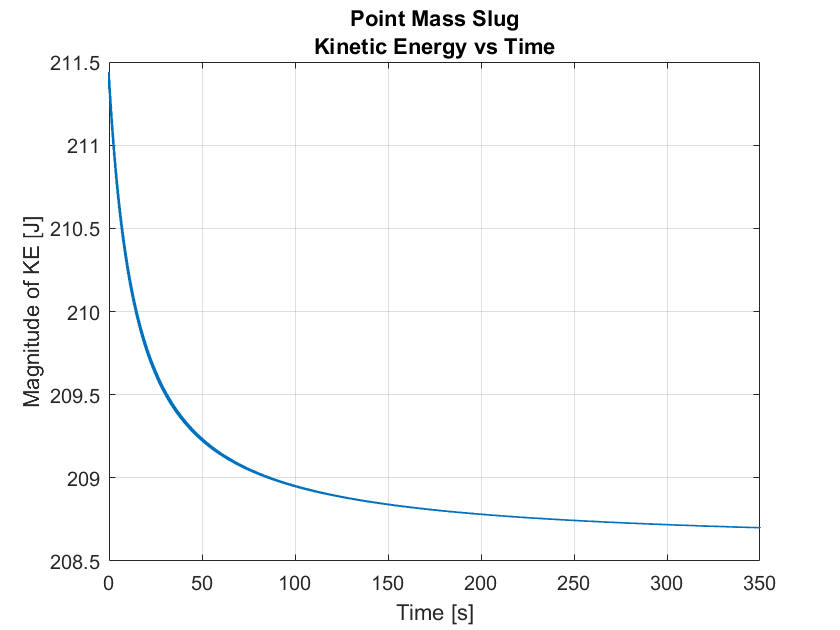}
\includegraphics[width=2in]{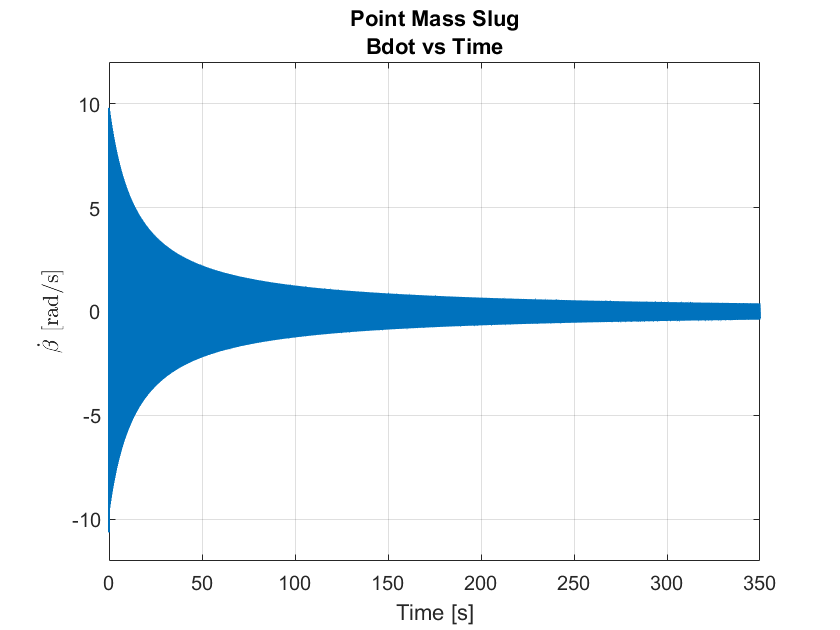}
\includegraphics[width=2.5in]{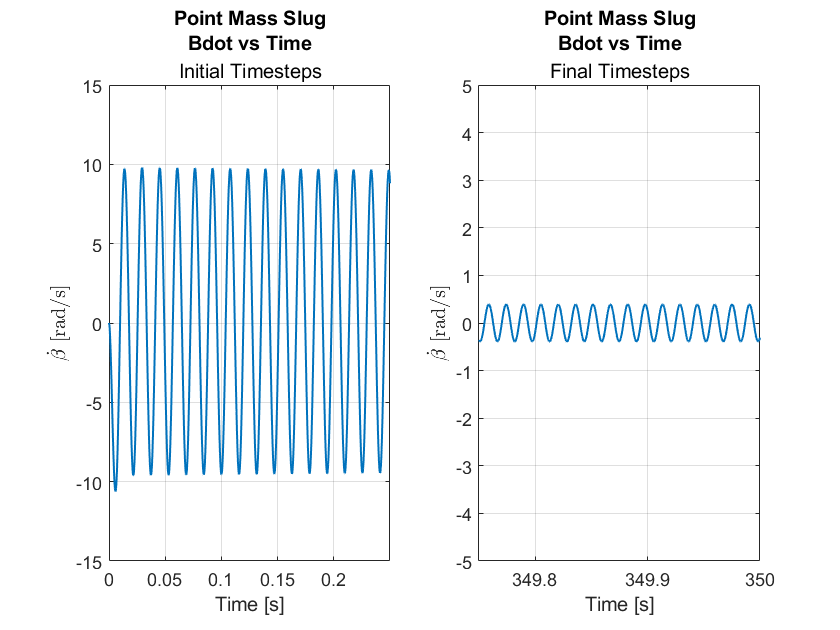}
\includegraphics[width=1.2in]{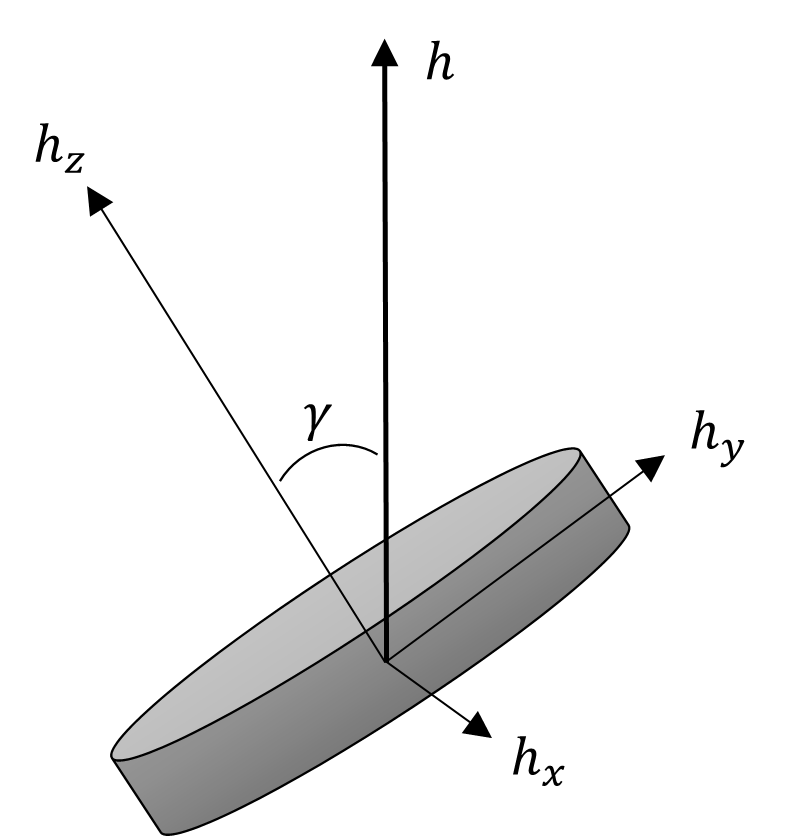}
\includegraphics[width=2.4in]{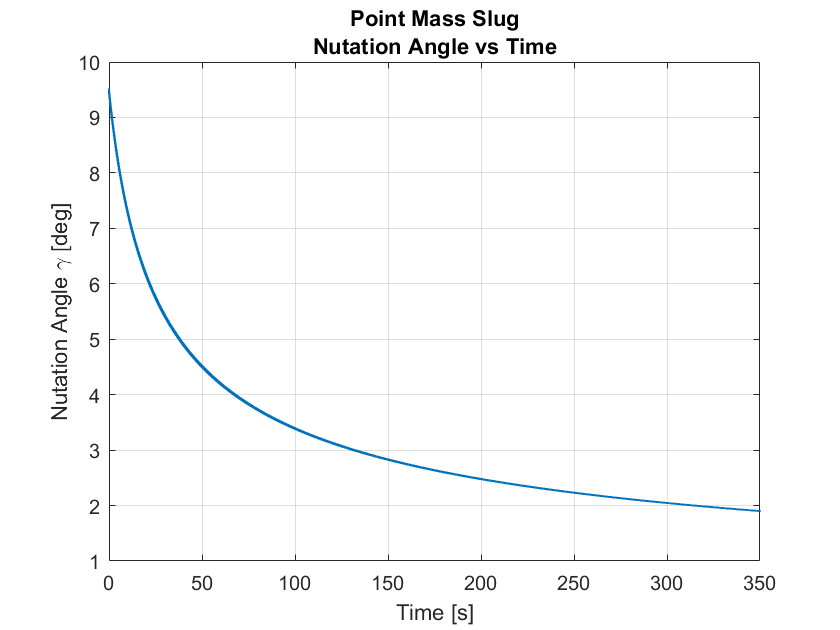}
\caption[Kinetic energy and nutation angle of the system in time (point-mass)]{Time plots showing the decay of various quantities in the point mass slug model: in the top row, kinetic energy and $\dot{\beta}$, on the original time interval. The bottom row left shows the $\dot{\beta}$ decay on two smaller time sub-intervals. The bottom right shows the decay of the nutation angle defined in the sketch.}
\label{betadotkenut}
\end{figure}
 In Fig. \ref{betadotkenut}, time plots of other quantities of interest are shown. The top left plot shows the decay of the system kinetic energy. The top right plot shows the behavior of the $\dot{\beta}$ variable on the original time interval and, in the bottom left plots, on two shorter subintervals. This variable represents the difference in angular velocity between the slug and the damper wall. The oscillating behavior seen suggests the slug transitions between rotational speeds faster and slower than the rotating body. This can be explained by the direction of the friction force acting on the slug. The direction of friction acts in the opposite the direction of $\dot{\beta}$ which would cause the velocity to approach zero. The size of the amplitudes corresponds to the amount of nutation present which would point to the transverse angular momentum accelerating the slug. This would cause oscillating behavior while friction acts against the motion. \\
\noindent{\bf Nutation damping.} In the absence of an external force such as gravity, the spacecraft body does not precess but can still experience nutation. Nutation is defined as the rotation of the spacecraft's major spin axis about the constant angular momentum vector $\bar{\mathbf{h}}$ in the spatially-fixed frame. Wlog, $\bar{\mathbf{h}}$ is assumed to be parallel to the $Z$-axis, as shown in Fig. \ref{rs}. The nutation angle $\gamma$ is defined as the angle between these two axes. Referred to the body-fixed frame, it is the angle between the $\mathbf{h}$ and $z$-axes, as shown in Fig. \ref{betadotkenut}, second row, and defined as $\tan(\gamma)=\frac{\sqrt{(h_x^2+h_y^2)}}{h_z}$. The nutation damper reduces the nutation angle such that spin vector approaches the desired spin axis (the $\bar{\mathbf{h}}$-axis) and removes the dynamic irregularities that arise from nutation.  
%
The decay of the nutation angle in the last plot of Fig. \ref{betadotkenut} shows that the system slowly restores the spin vector to desired spin axis. The point mass model does not fully dampen the nutation but is successful at significantly reducing nutation. The large settling time is undesirable as the system would take three hundred seconds to reach a nutation angle under two degrees. Parameters of the damper need to be tuned to maximize performance. \\
\noindent{\bf Vertical offset, radius and mass variations.} Literature suggests that offsetting the damper improves performance. The assumption that the center of mass is not displaced vertically can be used to add vertical offset $h_i$ as a parameter to analyze the performance. This would change the moment of inertia elements of the slug as follows $I_{xs} \rightarrow I_{xs} +m_rh_i^2, I_{ys} \rightarrow I_{ys} +m_rh_i^2$. However, no significant changes in damper performance were observed. 
On the other hand, it was observed that maximizing the mass and radius of the damper improves the settling time and reduces the approaching nutation angle. Given the low mass of CubeSats, however,  the choice of maximizing mass would not be desirable as the damper should minimize its effect on the dynamics of the satellite body. The values that performed the best within the range of parameters were found to be: Slug mass=0.01 kg and Damper Radius=0.05 m.

\section{Distributed slug.}
\begin{figure}[h]
\centering
\includegraphics[width=2in]{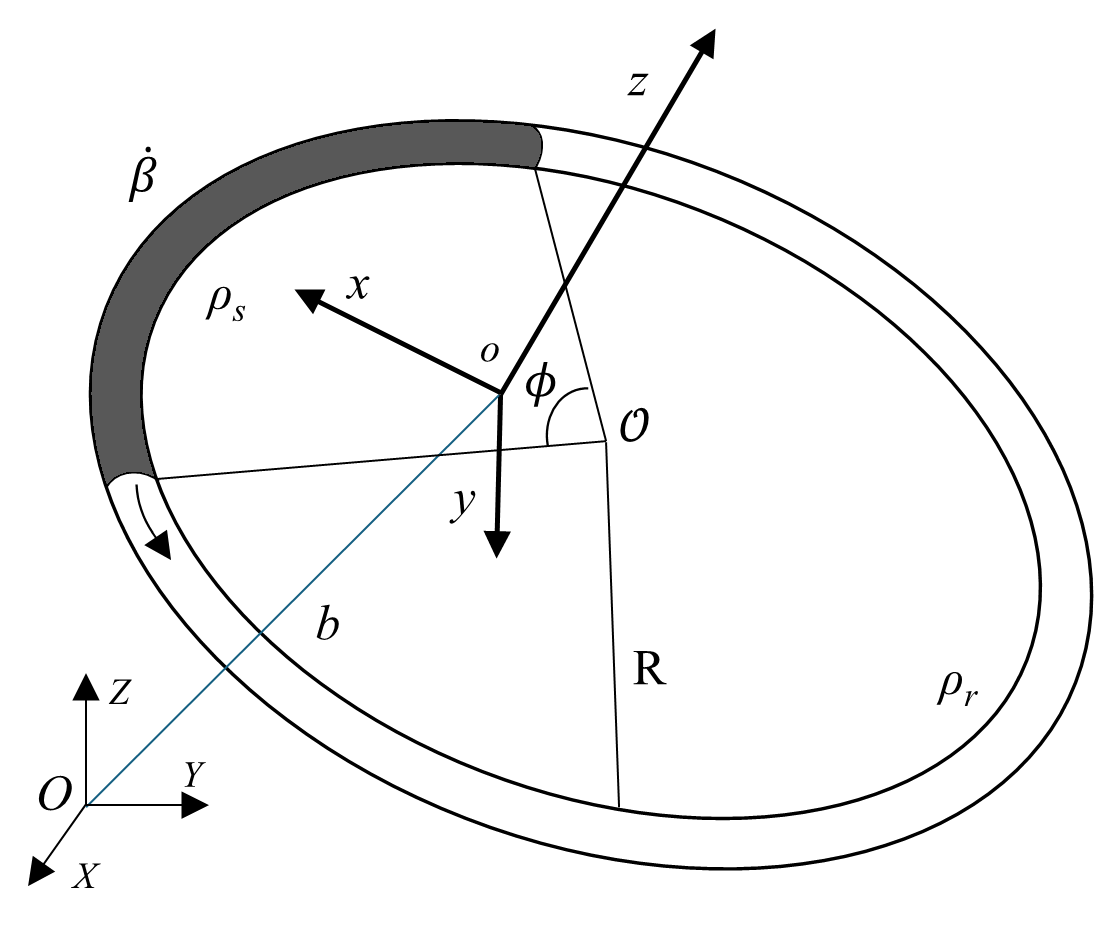}
\caption[Schematic sketch of the distributed ring damper]{Schematic sketch of the ring damper with fluid slug modeled as a distributed rigid mass.}
\label{rsdist}
\end{figure}

In this section, the slug is modeled more realistically as a distributed mass, as shown in Fig. \ref{rsdist}. But the rigid assumption is maintained so the slug is now a distributed rigid mass with no fluid flow. All other assumptions made in the point mass model continue to hold.  
\subsection{Dynamical Framework and Theoretical Analysis.}
The new geometry of the slug is (a portion of) a torus. The distribution of the slug is determined by the fill angle $\phi$, with $\phi=360$ degrees corresponding to a completely filled ring. Let the radius of the circular cross section of the ring be denoted by $r_c$ and let $R$ denote the radius of the ring.

The distributed mass leads to different inertia matrices for the slug.  The angular momentum of the slug in the body-fixed frame is now given by $
\mathbf{h}_s(t)=\mathbf{I_{ds}} \mathbf{\Omega_s}$,
where 
\begin{align*}
\mathbf{I_{ds}}=
\begin{bmatrix}
\rho_s \pi r_c^2  \frac{R^3}{2}(\phi -\sin\phi)  & 0 & 0\\
0 & \rho_s \pi r_c^2 R((d^2 + \frac{R^2}{2})\phi-4Rd \sin (\frac{\phi}{2})+\frac{R^2}{2} \sin (\phi)),& 0\\
0 & 0 & \rho_s \pi r_c^2 R ((d^2+R^2) \phi - 4Rd \sin (\frac{\phi}{2}) 
\end{bmatrix}
\end{align*}
(the `ds' stands for distributed slug).
The point mass slug inertia tensor (\ref{eq:pmslugI}) is recovered under the following limits: $\phi \rightarrow 0, \quad \sin \phi \rightarrow \phi, \quad \rho_s \pi r_c^2 R \phi \rightarrow m_s$. 
Finally, the expressions for the system angular momentum and kinetic energy in the body-fixed frame for the distributed rigid slug case are:
\begin{align*}
\mathbf{h}_{ds}&=\mathbf{I_{r}} \mathbf{\Omega}+\mathbf{I_{ds}} \mathbf{\Omega_s}, \quad
K.E._{ds}=\frac{1}{2}(\mathbf{\Omega})^T \mathbf{I_r} \mathbf{\Omega}+\frac{1}{2}(\mathbf{\Omega}_s)^T  \mathbf{I_{ds}} \mathbf{\Omega}_s
\end{align*}
\paragraph{Vertical offset.} With vertical offset in the distributed slug case, the following changes are seen in the inertia tensor. In particular, not all off-diagonal terms are equal to zero. This results in a few additional terms when the equations of motion are written in terms of their components, but the general form of the equations of motion remain the same. 
\begin{align*}
\bf{I_{ds}}&=\begin{bmatrix}
I_{xds} +m_rh_i^2 & 0 & I_{xz}\\
0 & I_{yds} +m_r h_i^2& 0\\
I_{xz} & 0 & I_{zds}
\end{bmatrix}, \quad I_{xz}=m_s h R_m \frac{\sin{(\frac{\phi}{2})}}{\frac{\phi}{2}}
\end{align*}
\paragraph{Equations of motion.}
The angular momentum equations are the same as (\ref{eq:angmomeqn2}) with $\mathbf{I_{ds}}$ replacing $\mathbf{I_s}$:
\begin{align}
\frac{d\mathbf{ h}}{dt}&= \mathbf{h} \times \left[(\mathbf{I_r}+\mathbf{I_{ds}})^{-1}  (\mathbf{h} +\mathbf{I_r}\dot{\beta} e_3)\right], \label{eq:angmomeqnds}
\end{align}
Following the same dissipation model as in the point mass slug, i.e. using equation (\ref{eq:omegazbeta}), results in
\begin{align*}
\ddot{\beta}&=-\dot{\beta}C_d (R-d)^2\frac{(I_{zs}+I_{zr})}{I_{zs}I_{zr}}+\frac{I_{xz}\dot{\Omega}_x}{I_{zs}}+\frac{d{h_z}}{dt}. 
\end{align*}
The terms involving $\dot{\Omega}_x$ and $dh_x/dt$ can be eliminated using $\dot{\mathbf{\Omega}}={(\bf{I_r}+\bf{I_s})}^{-1} (\dot{\mathbf{h}}-\mathbf{I_s} \ddot{\beta})$
and (\ref{eq:angmomeqnds}), leading to a final ODE in $\dot{\beta}$. 

%

\subsection{Projected Phase Portraits and Time Plots.}
The model no longer being a concentrated point mass required some changes to initial parameters. The damper radius was decreased to 0.025m to make sense of the fill angle and low mass of the fluid. The choice of a fill angle of 135 degrees was chosen from optimized parameters offered by previous literature \cite{alfriend1983comparison}. A vertical offset of $h_i=$ 0.015m was added to determine its effect on the performance of the fluid as previously mentioned. All other parameter values were chosen to be the same as in the point mass model and the trajectories were produced with the same disturbance and initial angular velocity. 
\begin{figure}[h]
\centering
\includegraphics[width=2.5in]{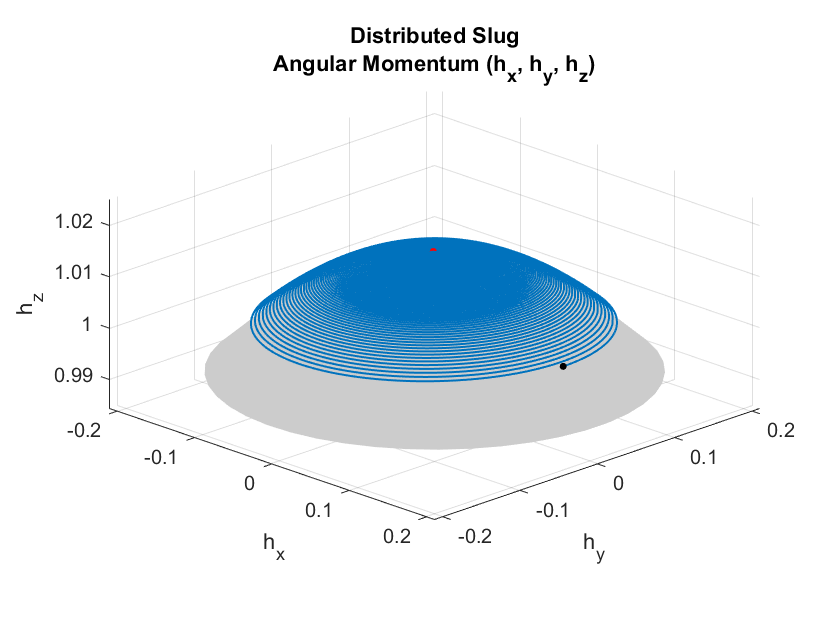}
\includegraphics[width=2.2in]{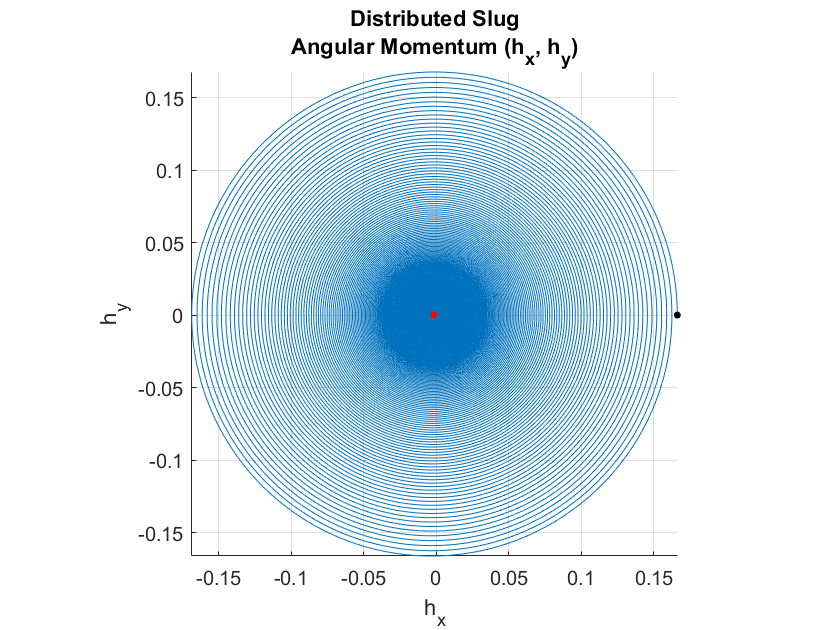}
\caption[Trajectory on the sphere of the distributed slug system]{A typical trajectory on the sphere of the dissipative system for the case of a distributed slug. Plot (B) being a 2-D projection on the $x-y$ plane showing the spiraling behavior of the system.}
\label{Spiral2}
\end{figure}

The trajectory on the sphere for the distributed slug displays an overall improvement over the damping of transverse angular momentum. In Fig. \ref{Spiral2}, the trajectory more clearly spirals towards the $z$-axis and does not show a significant nutation limit. The vertical offset now impacts performance as the off-diagonal inertia coefficients cause the slug to be affected by additional rotational accelerations. The simulation without a vertical offset displayed similar performance as the point mass model.  
\begin{figure}[h]
\centering
\includegraphics[width=2.4in]{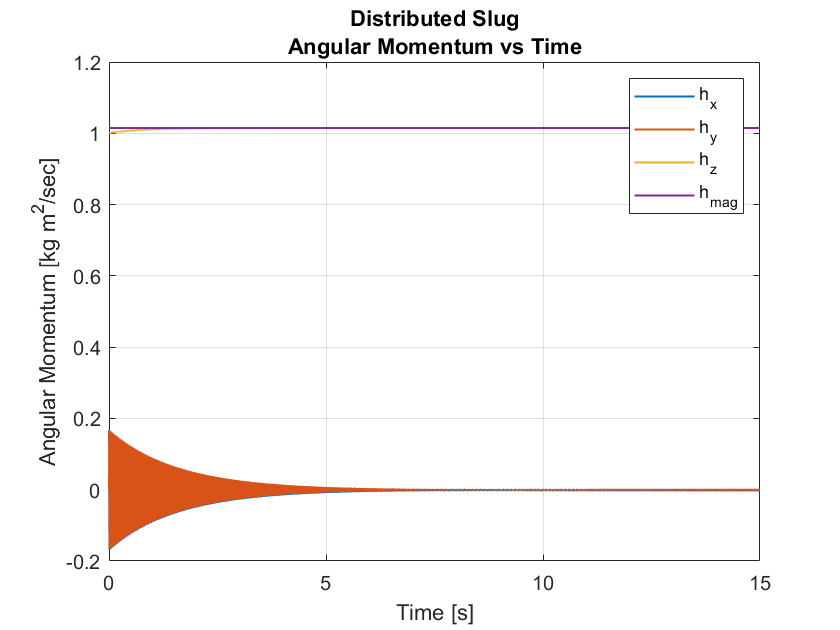}
\includegraphics[width=2.5in]{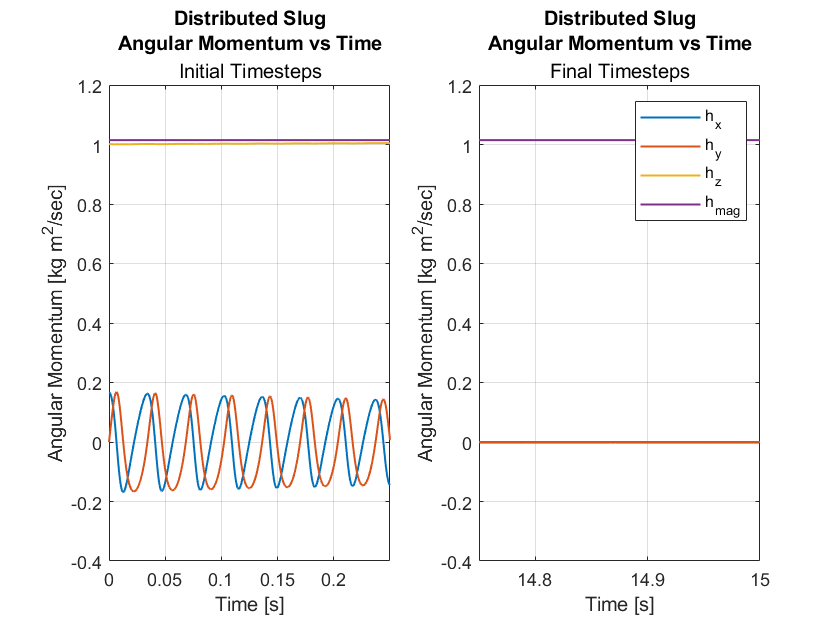}
\caption[Angular momenta components in time (distributed)]{Plots showing the conservation of the total magnitude of the angular momentum and the variation of each of its components in time. The almost full damping of the transverse angular momentum is clearly seen.}
\label{angmomds}
\end{figure}
The angular momentum vs time plots shown in Fig. \ref{angmomds} reflects the improvement in settling time over the point mass model. The transverse angular momentum approaches an almost zero value as opposed to the residual oscillations seen in the point mass model.  
\begin{figure}[h]
\centering
\includegraphics[width=2in]{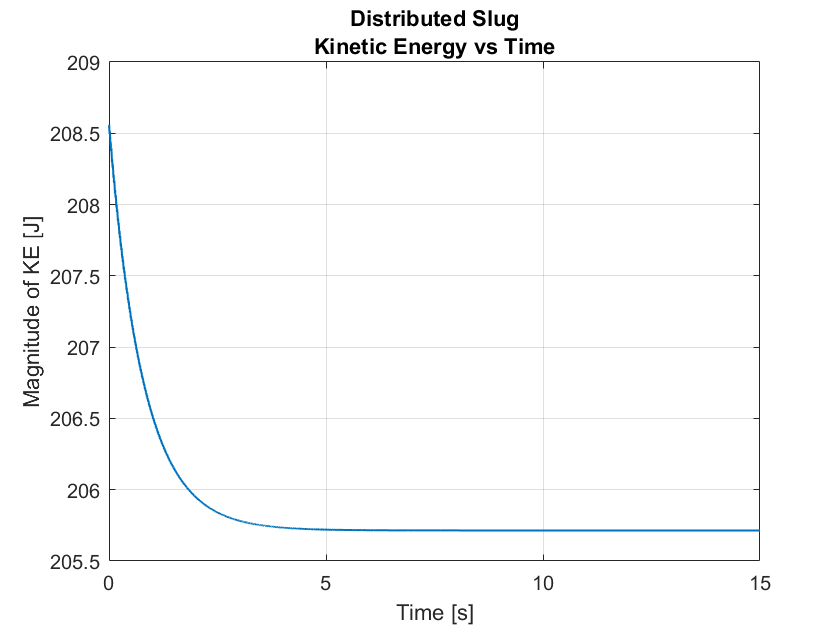}
\includegraphics[width=2in]{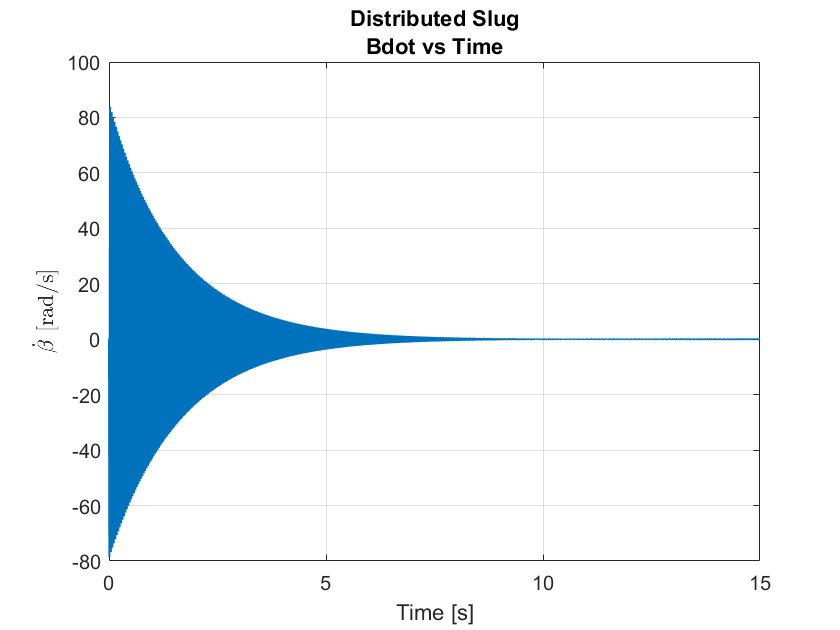}
\includegraphics[width=2.5in]{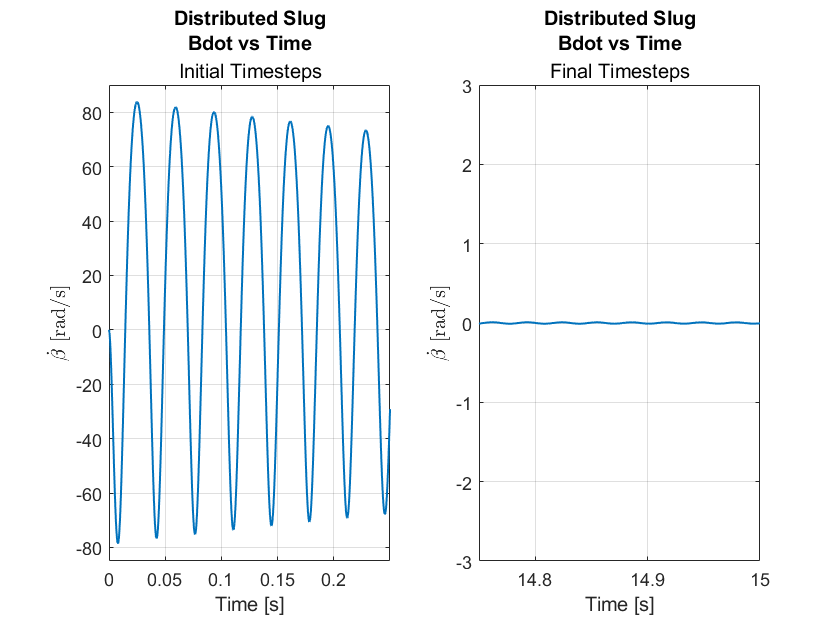}
\includegraphics[width=2.4in]{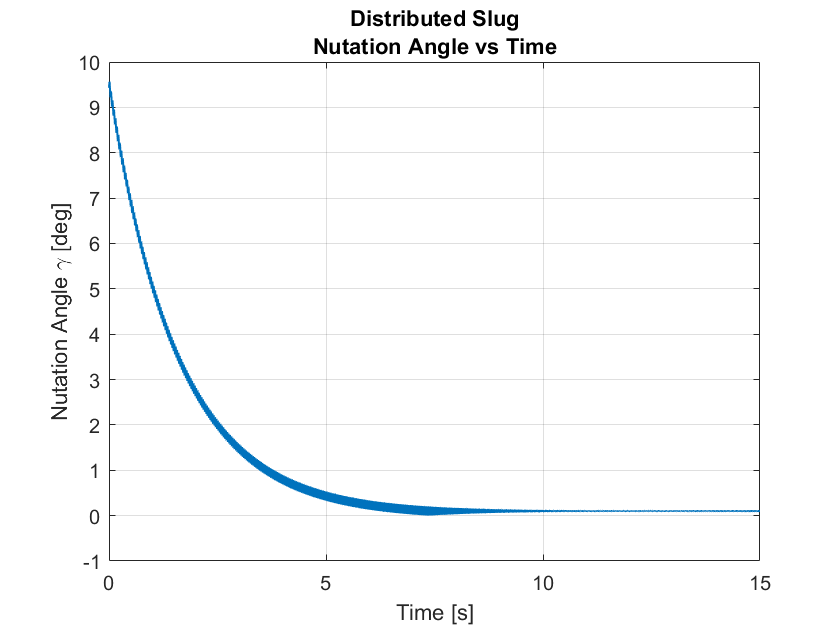}
\caption[Kinetic energy and nutation angle of the system in time (distributed)]{Time plots showing the decay of various quantities in the distributed slug model: in the top row, kinetic energy and $\dot{\beta}$, on the original time interval. The bottom row left shows the $\dot{\beta}$ decay on two smaller time sub-intervals. The bottom right shows the decay of the nutation angle. This shows the distributed damper's superior ability to remove the transverse angular momentum at a desired rate and have almost zero residual nutation.}
\end{figure}
The nutation angle decay reaches a steady point at about ten seconds after the disturbance. The small angle nutation damping mode can be seen in the slight oscillating behavior seen as the angle approaches zero. The low settling time now makes the possible integration to a small satellite relevant. 
%
For the point mass model, increasing parameters such as radius and mass showed a clear improvement in damping rate. For the distributed slug model, it was observed that the effect of increasing mass  on the approaching nutation angle was to decrease the settling time but increase the residual nutation angle. Parameters would need to be chosen optimally to minimize both  settling time and the nutation angle limit.
\section{CFD simulations of the flow of partially-filled slugs.}
\begin{figure}[h]
\centering
\includegraphics[width=2.3in]{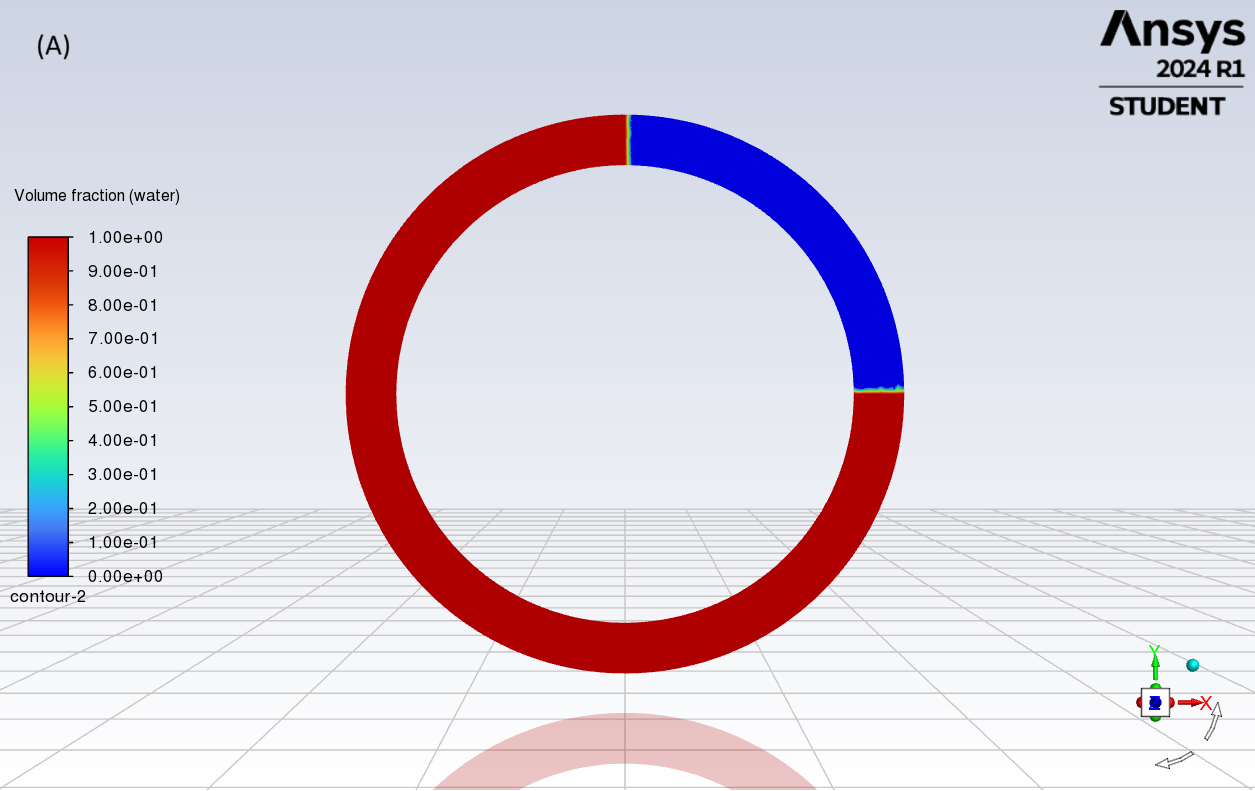}
\includegraphics[width=2in]{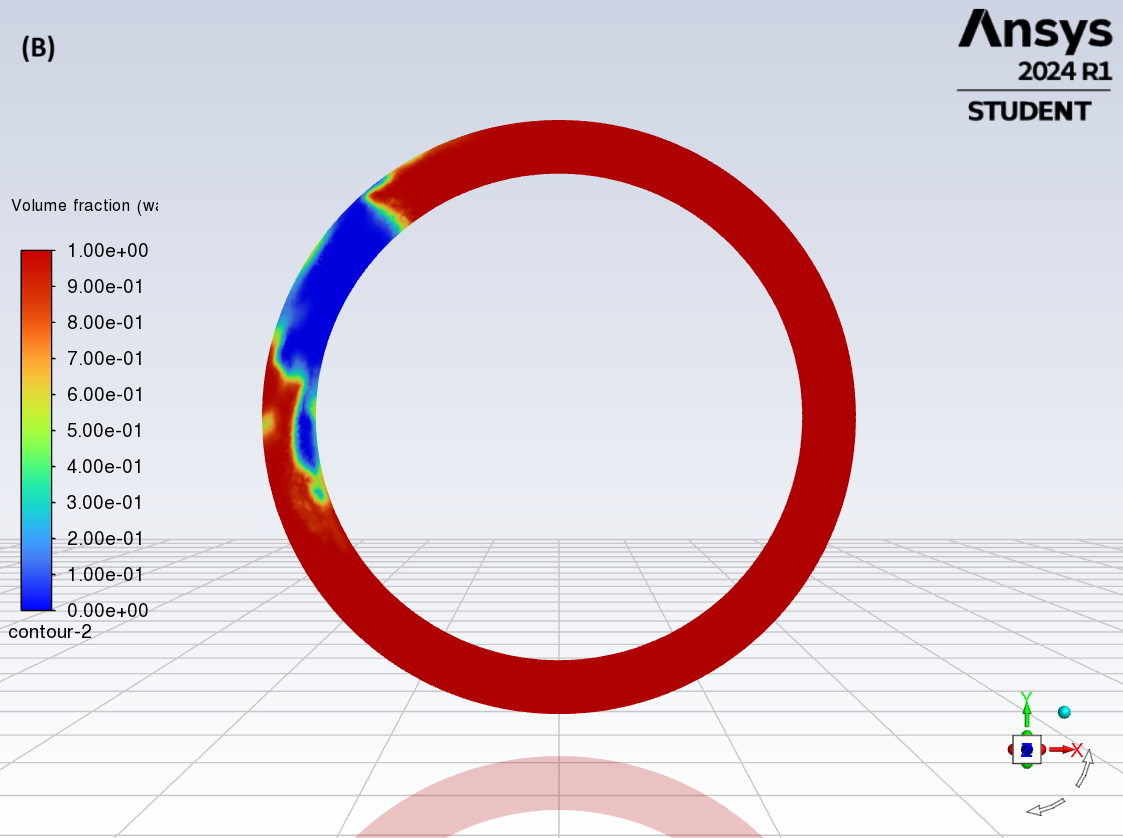}
\includegraphics[width=2in]{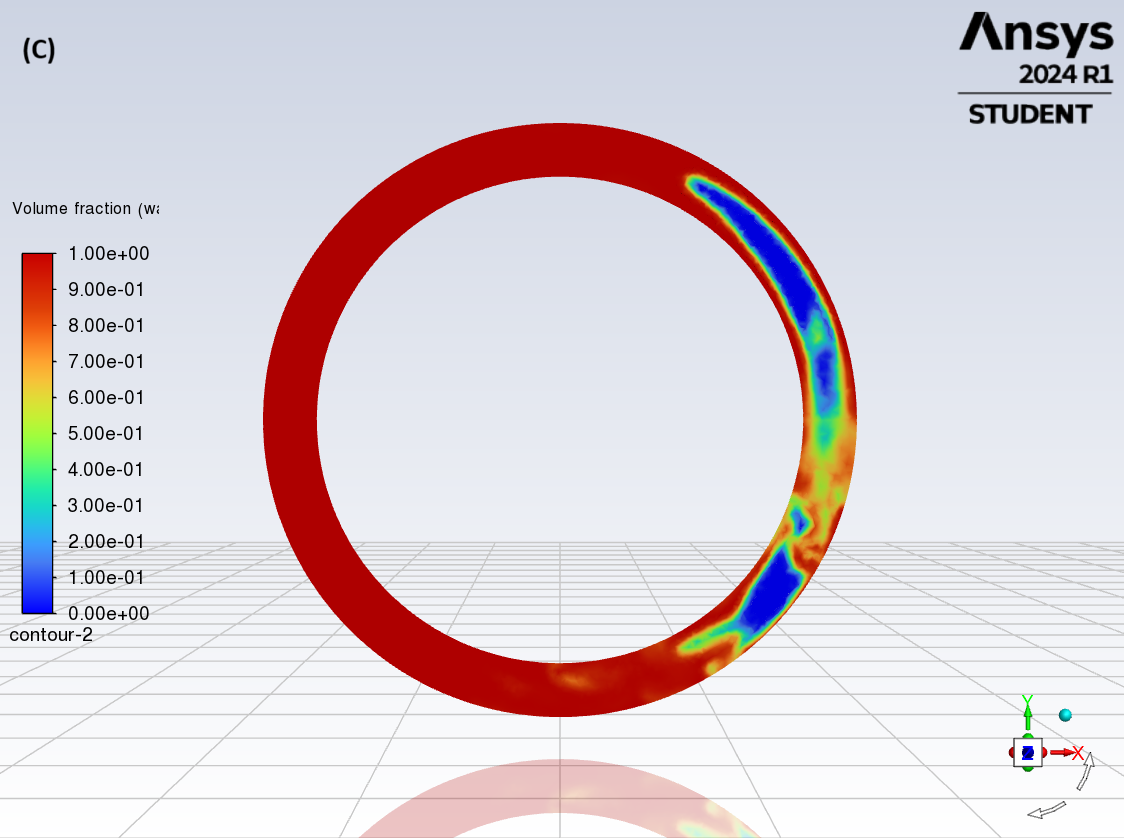}
\includegraphics[width=2in]{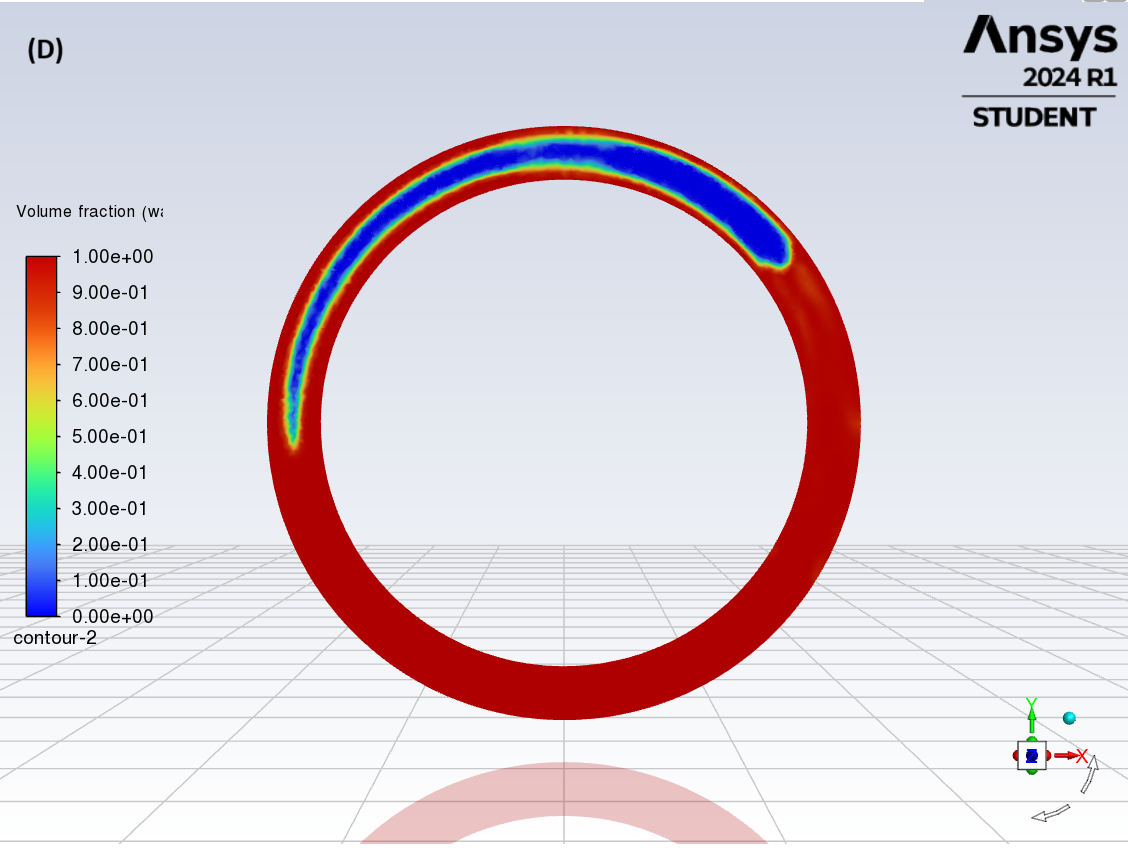}
\includegraphics[width=2in]{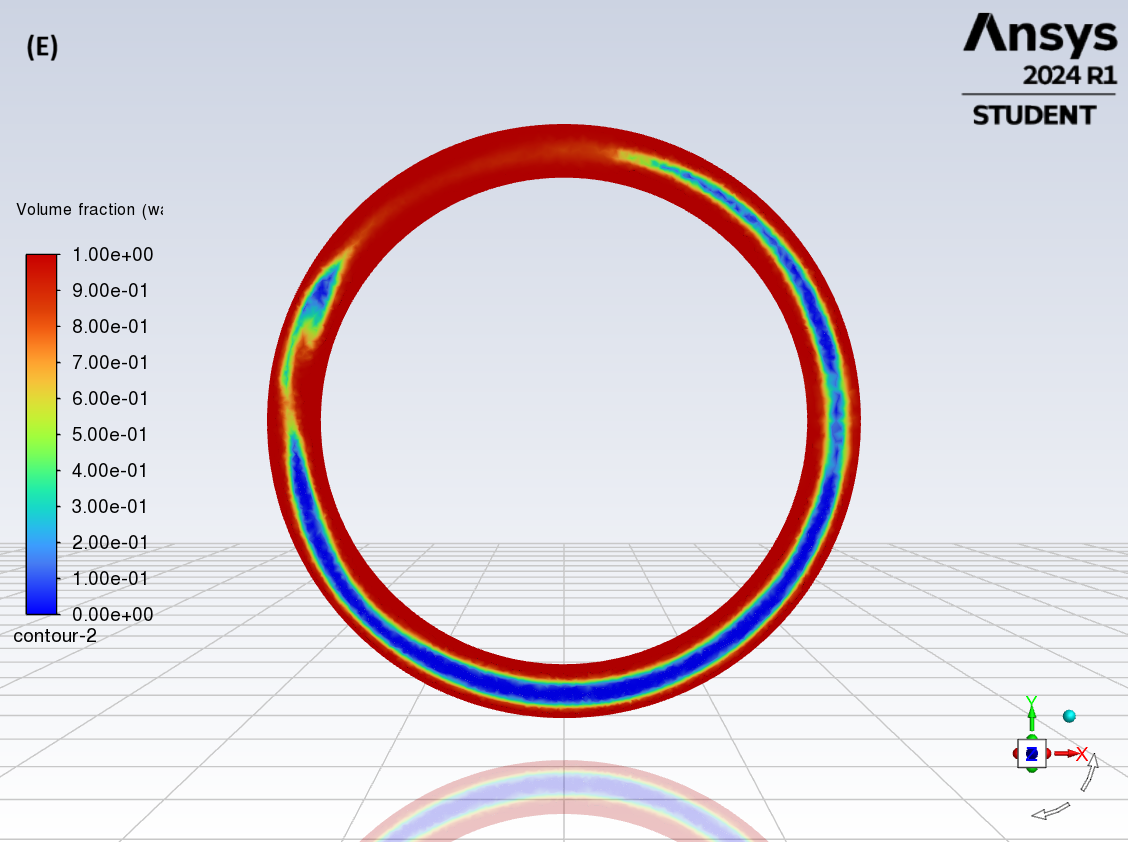}
\includegraphics[width=2in]{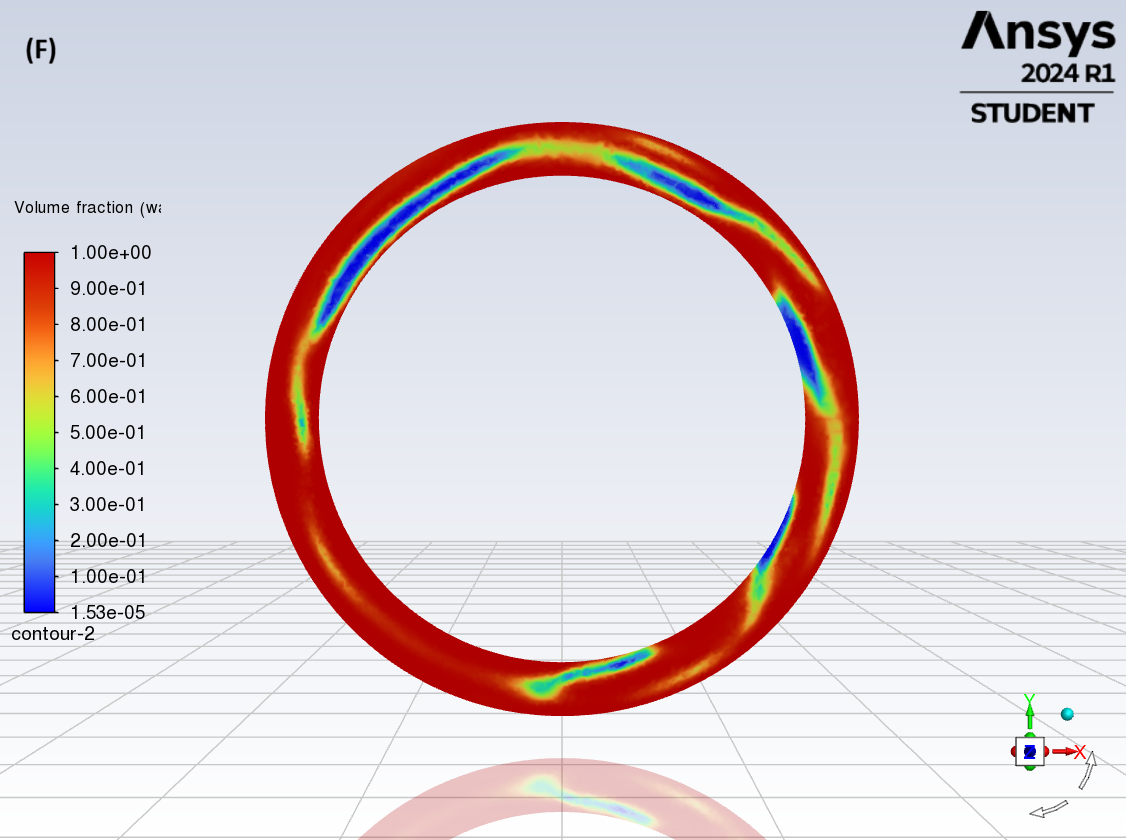}
\includegraphics[width=2.5in]{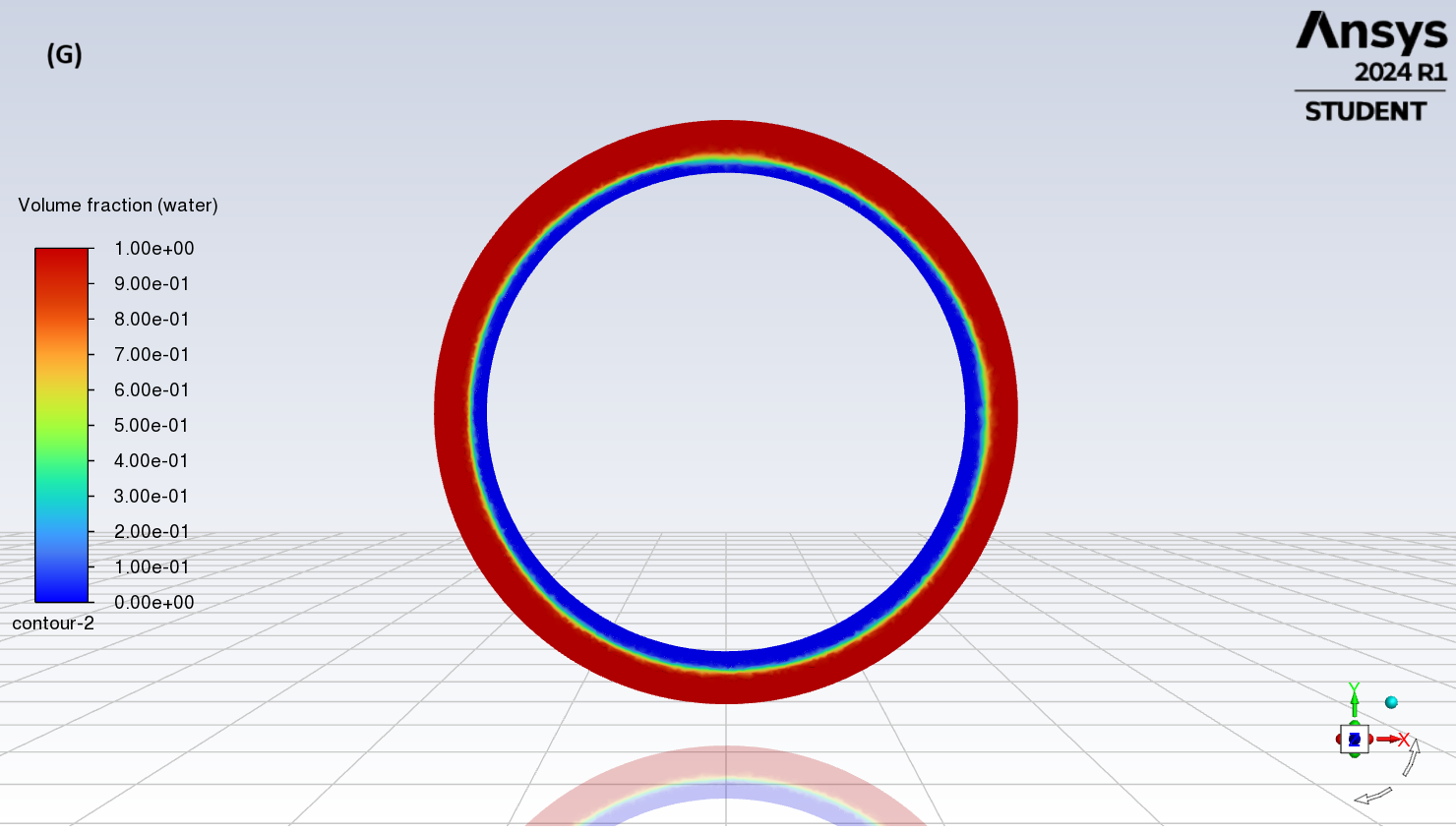}
\caption[Ansys Fluent: Simulation liquid and gas distribution timestamps]{The liquid (red) and gas (blue) distribution at various time instants of a CFD simulation of uniform ring rotation using the Ansys Fluent software.}
\label{snapshots}
\end{figure}
 CFD simulations were performed using the Ansys Fluent Fluid Simulation Software with the 2024 R1 Ansys Student version. The goal of the study is to obtain a better understanding of the complex fluid flow  generated within a nutation damper during rotation, with no simplifying assumptions made about the fluid slug. The geometry of the damper was modeled through the Ansys DesignModeler. 

    The primary sketch consists of a circular section of diameter 1 cm which represents the cross section of the damper. The revolve tool is used to create the annular shape where the geometry is created by revolving about the $z$-axis 5cm from the origin. The material chosen was Aluminum Alloy 2024 as it is commonly used in small satellites due to its relative low weight and strength. Portions of the ring were then separated to simulate a partially filled damper where the smaller cavity was then  filled with a gas. The main portion of the damper is filled with a liquid taking up a fill angle of 270 degrees. This was chosen as an average value from previous research on nutation dampers. Prior to the initialization the boundary conditions and material properties are set to create the conditions for the study. The ring body is given a moving wall boundary condition corresponding to a constant rotation rate of 50 rad/s about the $z$-axis. The speed was chosen to reflect high speed rotations seen in experimental research. The liquid and gas chosen were water and air, respectively. These materials were chosen because of their well-known properties and typical use. Liquids of high viscosity are typically used in dampers but in previous attempts the use of other liquids caused errors that did not allow for a full calculation. The first capture in Fig. \ref{snapshots} shows the distinct separate cavities that serve as the initial conditions for the study.
The resulting fluid distributions at various time instants of a CFD simulation can be seen in Fig. \ref{snapshots}. The heavier liquid gets pushed out to the outer wall of the damper which matches experimental observations \cite{cartwright1963circular}. The slug shape that is expected to be seen separates into distinct ring layers of liquid and gas. 
\begin{figure}[h]
\centering
\includegraphics[width=2.5in]{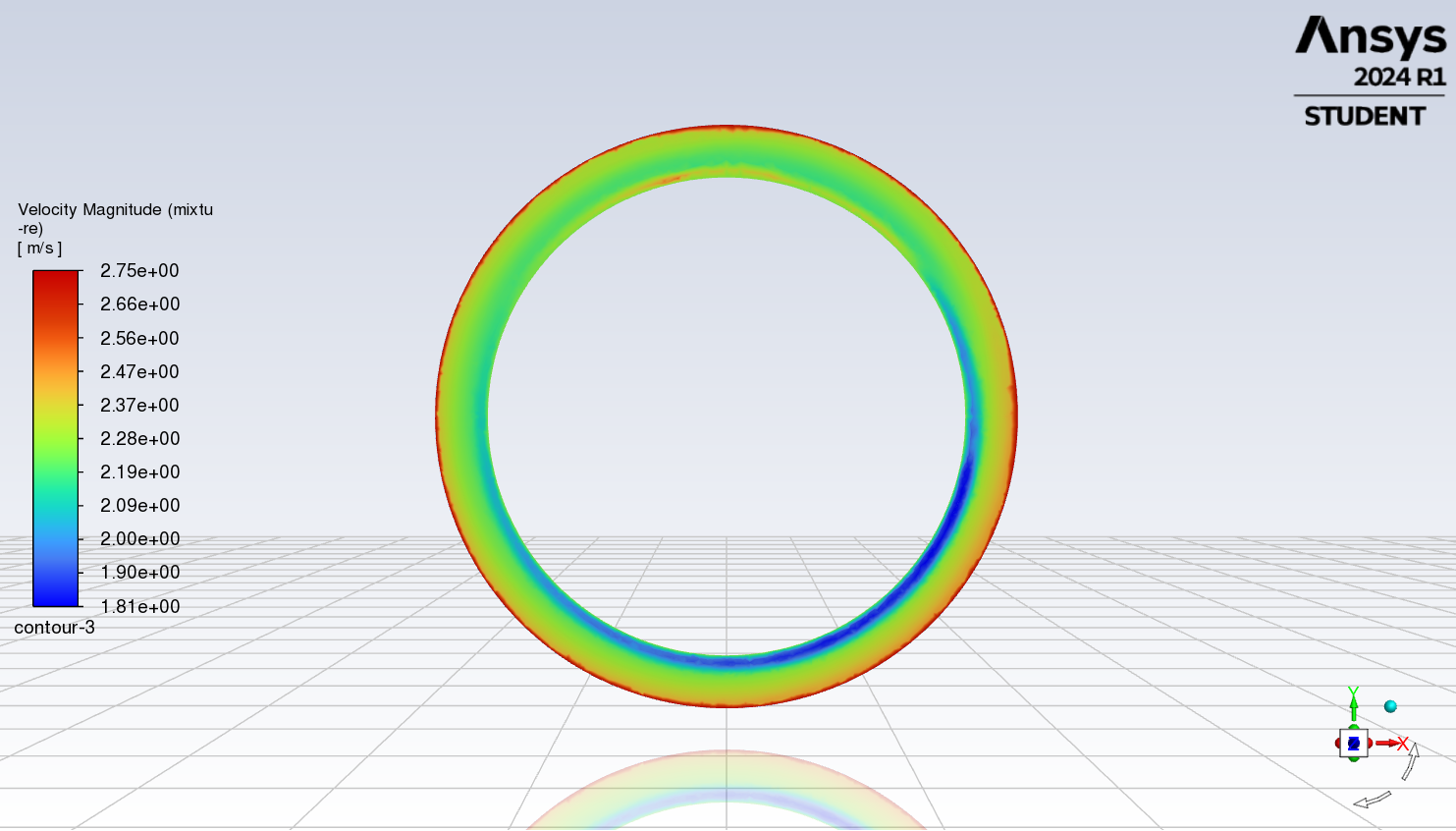}
\includegraphics[width=2.2in]{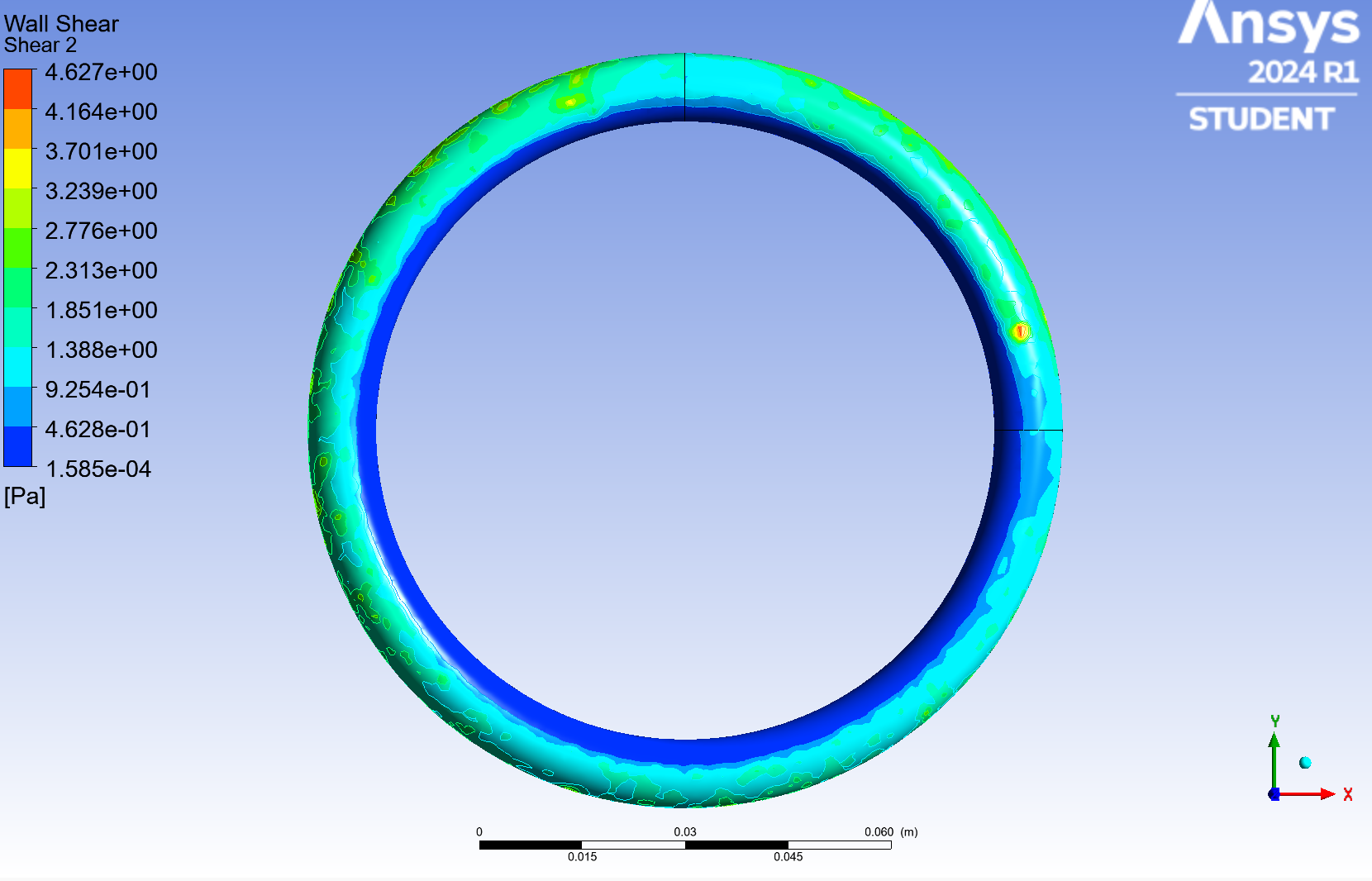}
\caption[Ansys Fluent: Velocity magnitude and wall shear stress distribution]{The velocity magnitude and the wall shear stress distribution at the end of the CFD simulation.}
\label{velshear}
\end{figure}
The distribution of the velocity magnitude in the fluid domain and the wall shear stress were obtained using the post solution analysis provided by the software, as shown in  Fig. \ref{velshear}. A surface integral over the surface area provides the skin friction coefficient in square meters. The resulting value of 0.0161m$^2$ divided by the total surface area provides the non-dimensional Cd value of 1.63. This was the value used in the MATLAB solutions as an estimate to the friction caused by the slug and the ring.

\section{Summary and future directions.}
 This paper examines the nonlinear dynamics of models for nutation dampers which are devices used in the satellite industry to dampen the unnecessary nutation wobbles induced in a spinning satellite upon encountering a disturbance. The geometry of the damper model is that of a rigid circular ring  partially filled by a liquid slug. Though these devices and models have been around for decades, the scientific literature on them is somewhat sparse. The problem is interesting from both theoretical and numerical perspectives. 
The purpose of this paper is twofold. It highlights the importance of the angular momentum sphere in understanding the dynamics. Secondly, it  makes use of modern computational power--not available to early researchers on the problem--to investigate the complex unsteady fluid flow in the ring. Though these numerical simulations are only preliminary, it may be considered a novel feature of this paper.  Phase portrait simulations and time histories of important parameters--such as the system kinetic energy and the nutation angle--and CFD simulations using the Ansys Fluent software are presented. 

   The basic research presented in this paper could be extended in several directions, both theoretical and with applications in mind. Focusing on engineering applications, one obvious direction is to explore more thoroughly applications to CubeSats and nanosatellites. There are many publications of proposed or ongoing CubeSat missions that utilize spin stabilization [22-24]. Another direction is to develop a Navier-Stokes CFD code that is coupled to the rigid body equations. This would allow one to investigate the fully dynamically coupled fluid-structure interaction problem with no simplifying assumptions made on the fluid flow in the slug, which seems to be a largely unexplored topic; in particular, to explore the effects of the non-uniform rotation that occurs in a real damper. Though a challenging task, such a code could be useful for applications in other areas of engineering as well, for example, the sloshing of liquids in partially-filled tubes.  

\listoffigures
\end{document}